\documentclass[12pt, reqno, a4paper]{amsart}
\usepackage{graphics}
\usepackage{subfigure}
\usepackage{amsmath,amssymb,amsthm}
\usepackage{enumerate}

\usepackage{prettyref,pretty_styles, theorem_styles}
\usepackage{general, powerProduct, GBAS}

\newcommand{\rtr}{\operatorname*{rtr}}

\newcommand{\restRel}[2]{{#1}_{\vert_{{#2}}}}
  
\newcommand{\Integ}{\Xi}
\newcommand{\Deriv}{\xi}

\newcommand{\shrink}{\operatorname*{int}}
\newcommand{\E}{\mathcal{E}}
\newcommand{\dual}[1]{{#1}^*}
\newcommand{\nicemonoid}{cancellative, torsion free, reduced abelian
  monoid}  
\renewcommand{\maxsupp}{\gamma}
\newcommand{\wrt}{wrt\ }

\begin{document}
\title[Partial orders associated to generic initial ideals]{On some partial
  orders associated to generic initial ideals} 
\author{Jan Snellman}

\begin{abstract}
  We study two partial orders on \([x_1,\dots,x_n]\), the free abelian
  monoid on \(\set{x_1,\dots,x_n}\). These partial orders, which we
  call the  ``strongly stable'' and the ``stable'' partial order,  are
  defined by the property that their filters are precisely the
  strongly stable and the stable monoid ideals. These ideals arise in
  the study of generic initial ideals.
\end{abstract}
\date{\today}
\address{Laboratorie GAGE\\
Ecole Polytechnique\\
91128 Palaiseau Cedex\\ 
France}
\email{jans@matematik.su.se}
\thanks{The author was supported by a grant from \textit{Svenska institutet}
  and by grant n. 231801F from \textit{Centre International des Etudiants et
  Stagiaires}.} 
\maketitle
\sloppy

\begin{section}{Introduction}
     So called \emph{strongly stable} (or \emph{Borel}) monomial ideals 
   are of interest because they appear as \emph{generic initial ideals}
   \cite{Galligo:Prep,Galligo:DivStab,Green:gin, Ebud:View}. 
   They have been used to give new proofs of the
   Macaulay and Gotzmann theorems for the growth of Hilbert series, and
   to extend these results to related rings
   \cite{Bigatti95, Aramova:Gotzman, Gasharov:Gotzmann}. A related
   class of ideals, the so called \emph{stable} monomial ideals, are
   also of interest \cite{Eliahou:MinRes}.

   A subset \(V \subset [x_1,\dots,x_n]_d\) (of the set of monomials of total
   degree \(d\) in \(n\) variables) 
   is \emph{strongly stable} (or \emph{Borel}) 
   if whenever a monomial \(m\) belongs to \(V\), then
   \(\frac{x_i}{x_j}m \in V\) for all \(1 \le i < j \le n\) such that
   \(\divides{x_j}{m}\).
   A monoid ideal \(I \subset [x_1,\dots,x_n]\) in the free abelian
   monoid on \(n\) letters is strongly stable iff
   \(I_d\) is strongly stable for all \(d\)). 
   The anti-symmetric binary relation on
   \([x_1,\dots,x_n]\) which consists of all such pairs 
   \((\frac{x_i}{x_j}m,m)\) is usually called the relation of
   \emph{elementary moves}. As observed in
   \cite{Marinari:SomeProp}, \(V\)
   is strongly stable iff it is a
   \emph{filter} \wrt the poset \(A_{n,d}\) which is
   the reflexive and  transitive closure of the elementary moves relation. 

   Also in \cite{Marinari:SomeProp}, the problem of determining all
   Borel monomial ideals with the same \(h\)-vector as
   of certain Artinian algebras is considered. The authors show that 
   the problem can be reduced to the enumeration
   of all Borel subsets of a fixed cardinality. Since Borel subsets are
   precisely the filters of \(A_{n,d}\), one is
   interested in studying the poset of filters of \(A_{n,d}\). We
   show that this latter poset is isomorphic to the poset of certain
   hyper-partitions. In particular, for \(n=3\), we obtain a bijection
   between Borel subsets of cardinality \(v\) of monomials of degree
   \(d\) and numerical partitions of \(v\) into distinct parts not
   exceeding \(d+1\).

   Before arriving at this result, we must conduct a study of the
   poset \(A_{n,d}\) itself.
   We find that \(A_{n,d}\) is a 
   distributive lattice, namely  the lattice
   of Ferrer's diagrams that fit into a \((n-1) \times d\) box. 
   It is therefore self-dual, ranked with the rank function given by
   the \(q\)-binomial coefficients, rank-symmetric, et cetera. 

   Next, we set out to make a partial order on \([x_1,\dots,x_n]\) such
   that its filters are precisely the Borel monoid ideals. To
   accomplish this, one must involve the 
   divisibility partial order \(D\), since a monoid ideal in
   \([x_1,\dots,x_n]\) is nothing but a filter \wrt
   \(D\). We denote by \(A_{n,\cdot}\) the reflexive-transitive
   closure of the  union of \(D\) and \(A_{n,d}\), for all \(d\). 
   Since it has been obtained by ``gluing'' together all
   \(A_{n,d}\) by divisibility, this poset has exactly the Borel
   monoid ideals as its filters. It is a distributive
   lattice intimately related to the Young lattice. Somewhat
   surprisingly, it turns out  that this poset is also the intersection
   of all term orders on  \([x_1,\dots,x_n]\) which restrict to the
   correct ordering of the variables. This result can be
   regarded as an marginal note to the well-known classification of
   term-orders \cite{TermOrderings}.
  
   Clearly, \(A_{n,\cdot}\) is in a
   natural way included in \(A_{n+1,\cdot}\). Passing to the
   inductive limit, we get yet another interesting poset, once again a
   distributive lattice, but this time not quite the same thing as the
   Young lattice; however, if we order our variables so that \(x_1\) is
   the \emph{smallest} rather than the largest variable, and carry out
   the above construction, we do get the full Young lattice.

   Similarly, one may consider the stable relation \(B_{n,d}\) on
   \([x_1,\dots,x_n]_d\). Here, only the
   smallest occurring variable in the monomial is allowed to be replaced
   with something larger. This seemingly inconsequential change yields
   a poset drastically different from the strongly stable poset: it is
   a lattice, which however is not distributive, not
   even modular, and which is probably not isomorphic to any of the
   classic posets. We give an explicit description of the meet
   operation on \(B_{n,d}\).

   As for the poset of filters of \(B_{n,d}\),
   we do not obtain such a nice description as for the strongly stable
   case, but for the special case \(n=3\), \(d > v\) we
   show that the number of stable subsets of cardinality \(v\) is equal to the
   number of 
   fountains of \(v\) coins.

   \begin{subsection}{Acknowledgement}
     The connection between \(A_{n,d}\) and the distributive
     lattice of Ferrer's diagrams that fit inside an \((n-1) \times d\)
     box was pointed out by Anders Bj\"orner.
   \end{subsection}

\end{section}

\begin{section}{Multiplicative relations and multiplicative partial
    orders}
  We denote by \(\Nat^+ = \set{1,2,\dots}\) the set of positive integers,
and by \(\Nat=\set{0,1,2,\dots}\) the set of non-negative integers.
  For relations, partial orders, lattices, and commutative semigroups,
  we will endeavour to follow the terminology in \cite{Gilmer:Semi,
    Birkhoff:LT, GLT}. 

  We mean by a binary relation \(R\) on a set \(M\) a subset of \(M \times
  M\). If \((a,b) \in R\) we sometimes write \(a R b\), or say that \(a \ge
  b\) with respect to \(R\).

  \begin{definition}
    For any binary relation \(R\), we denote by \(\rtr(R)\) the
    reflexive and transitive closure of \(R\).
  \end{definition}

  \begin{definition}
    If \(R\) is a binary relation on a set \(A\), and if \(B \subset
    A\), then we denote the restriction of \(R\) to \(B\) by
    \(\restRel{R}{B}\).
  \end{definition}

  \begin{definition}
    If \(M\) is a \nicemonoid,
    and \(R\) is an anti-symmetric binary relation on \(M\), then
    \(R\) is said to be \emph{(strongly) multiplicative} if the
    following condition holds:
    \begin{equation}
      \label{eqn:strongmult}
      \forall m, m', t \in M: \quad (m, m') \in R \quad \implies \quad
      (mt, m' t) \in R
    \end{equation}
    A relation \(S\) on \(M\) which is contained in some strongly
    multiplicative 
    relation is called \emph{weakly multiplicative}.
  \end{definition}

  We will need the following theorem, which follows immediately from
  \cite[Corollary 3.5]{Gilmer:Semi}
  \begin{thm}\label{thm:gilmerplus}
    Let \(M\) be a abelian, cancellative, torsion-free monoid, and let \(R\)
    be a strongly multiplicative partial order on \(M\). Suppose that \(x,y
    \in M\) and that
    \(\set{x^n,y^n}\) is an antichain in \(M\) for all \(n \in \Nat^+\). Then
    \(R\) can be extended to a 
    multiplicative total order \(\tilde{R}\) such that \((x,y) \in
    \tilde{R}\). 
  \end{thm}

\end{section}

\begin{section}{Definition of the stable and the strongly stable
    partial order}
    Let \(X=\set{x_1,x_2,x_3,\dots}\), and let, for \(n \in \Nat^+\), 
  \(X^n =  \set{x_1,\dots,x_n}\). Denote by \(\M\) the free abelian monoid on
  \(X\), and by \(\M^n\) or \([x_1,\dots,x_n]\) the free abelian
  monoid on \(X^n\). 
  For \(d  \in \Nat\) we denote by \(\M_d\) and \(\M_d^n\) the subsets
  consisting of elements of total degree \(d\). We will occasionally
  use \([x_1,\dots,x_n]_d\) as a synonym for \(\M_d^n\). If \(m \in
  \M\) then 
  we define \(\supp(m) = \setsuchas{i \in \Nat^+}{\divides{x_i}{m}}\),
  and \(\maxsupp(m) = \max( \supp(m))\), with the convention that
  \(\maxsupp(1)=0\).  
  
  \begin{definition}\label{def:strongstab}
    Let \(n,d \in \Nat\). The \emph{strongly stable partial order}
    \(A_{n,d}\) on
    \(\M_d^n\) is the reflexive-transitive closure of the set of all pairs 
    \begin{equation}
      \label{eqn:strongstab}
      A_{n,d}^\circ = \setsuchas{(m,m') \in \M_d^n \times \M_d^n}{\exists
      i,j: i < j , \,\, m = \frac{x_i}{x_j} m'}
    \end{equation}
  \end{definition}

  \begin{definition}\label{def:stab}
    Let \(n,d \in \Nat\). The \emph{stable partial order} \(B_{n,d}\) on
    \(\M_d^n\) is the reflexive-transitive closure of the set of all pairs  
    \begin{equation}
      \label{eqn:stab}
      B_{n,d}^\circ = \setsuchas{(m,m') \in \M_d^n \times \M_d^n}{\exists i:
      \, i
      < \maxsupp(m'), \,\, m = \frac{x_i}{x_{\maxsupp(m')}} m'}
    \end{equation}
  \end{definition}
  
  So \(x_1^2 x_3 \ge x_1 x_2 x_3\) with respect to \(A_{3,3}^\circ\), but not
  with respect to \(B_{3,3}^\circ\).
  
  \begin{lemma}
    \(A_{n,d}\) and \(B_{n,d}\) are partial orders. \(B_{n,d} \subset
    A_{n,d}\), and the inclusion is strict for \(n>2\), \(d>2\).
  \end{lemma}
  \begin{proof}
    The last assertions are obvious, and hence it is enough to show
    that \(A_{n,d}\) is anti-symmetric. Let \(f: \M_d^n \to \Nat\) be
    defined by \[f(x_1^{\alpha_1} \cdots x_n^{\alpha_n}) = \alpha_1 +
    2\alpha_2 + \cdots + n \alpha_n.\] Then every substitution
    \[x_1^{\alpha_1} \cdots x_n^{\alpha_n} \to x_1^{\alpha_1} \cdots
    x_{i-1}^{1+\alpha_{i-1}} x_{i}^{\alpha_i - 1}
    \cdots x_n^{\alpha_n}\] strictly lowers the \(f\)-value. Hence, 
    the reflexive-transitive closure of \(A_{n,d}^\circ\) is anti-symmetric.
  \end{proof}

  We observe that \(x_1^d\) and \(x_n^d\) are the unique maximal and minimal
  elements of both \(A_{n,d}\) and \(B_{n,d}\).

  \begin{definition}\label{def:D}
    We denote the divisibility partial order on \(\M\) by
    \(D\). Abusing our notations, we denote any restriction of \(D\)
    to a subset of \(\M\) simply by \(D\). 
  \end{definition}
  
  Thus \(m \le m'\) with respect to \(D\) if and only if \(\divides{m}{m'}\),
  that is, if and only if \(m' = tm\) for some \(t\).
  
  \begin{prop}\label{prop:inclu}
    For \(d,n,n' \in \Nat^+\), \(n \le n'\),  we have that
    \(\restRel{A_{n',d}}{\M_d^n} = A_{n,d}\), and similarly for
    \(B\). 
  \end{prop}
  
  \begin{theorem}\label{thm:glue}
    The sets
    \begin{align}      \label{eqn:glueAB}
      A_{\cdot,\cdot} & := \rtr(D \cup \bigcup_{n,d} A_{n,d}) \\
      B_{\cdot,\cdot} & := \rtr(D \cup \bigcup_{n,d} B_{n,d}) 
    \end{align}
    are partial orders on \(\M\); \(A_{\cdot,\cdot}\) is
    strongly multiplicative and contains \(B_{\cdot,\cdot}\), which is
    therefore weakly multiplicative. We define the following restrictions:
    \begin{align}\label{eqn:mglueAB}
      A_{\cdot,d} &= \restRel{A_{\cdot,\cdot}}{\M_d} = \bigcup_{n}
      A_{n,d} \subset \M_d \times \M_d\\
      B_{\cdot,d} &= \restRel{B_{\cdot,\cdot}}{\M_d} = \bigcup_{n}
      B_{n,d} \subset \M_d \times \M_d \\
      A_{n,\cdot} &= \restRel{A_{\cdot,\cdot}}{\M^n} =
      D \cup \bigcup_{d} A_{n,d} \subset \M^n \times \M^n \\
      B_{n,\cdot} &= \restRel{B_{\cdot,\cdot}}{\M^n} =
      D \cup \bigcup_{d} B_{n,d} \subset \M^n \times \M^n
    \end{align}
    Then \(A_{n,\cdot}\) is strongly multiplicative, whereas
    \(B_{n,\cdot}\) is weakly multiplicative.
  \end{theorem}
  \begin{proof}
    It is enough to show that \(A_{\cdot, \cdot}\) and \(B_{\cdot,
    \cdot}\) are partial orders, and that \(A_{\cdot, \cdot}\) is
    strongly multiplicative. Define a function 
    \begin{align*}
    g: \M &\to \Z \times \Z\\
    g(x_1^{\alpha_1} \cdots x_n^{\alpha_n}) &= 
    (-\alpha_1-\cdots-\alpha_n, \alpha_1 +
    2\alpha_2 + \cdots + n \alpha_n),      
    \end{align*}
    and order \(\Z \times \Z\)
    lexicographically.
    Then substitutions of the form
    \begin{align*}
    x_1^{\alpha_1} \cdots x_n^{\alpha_n} &\mapsto x_1^{\alpha_1} \cdots
    x_{i-1}^{1+\alpha_{i-1}} x_{i}^{\alpha_i - 1}
    \cdots x_n^{\alpha_n}\\ 
    x_1^{\alpha_1} \cdots x_i^{\alpha_i} \cdots x_n^{\alpha_n} &\mapsto 
    x_1^{\alpha_1} \cdots x_i^{1+\alpha_i} \cdots x_n^{\alpha_n}      
    \end{align*}
    leads to a \(g\)-value which is strictly smaller. Hence,
    \(A_{\cdot, \cdot}\) is anti-symmetric and therefore a partial
    order.

    If \((m,m') \in A_{\cdot, \cdot}\) then there is a finite chain
    \begin{displaymath}
      m \to m_1 \to m_2 \to \dots \to m'
    \end{displaymath}
    where each arrow is a substitution of one of the two types
    above. It is easily seen that both such substitutions remain valid
    after multiplication by an arbitrary element \(t \in \M\), hence
    \((tm,tm') \in A_{\cdot, \cdot}\). 
  \end{proof}

\end{section}

\begin{section}{The r\^aison d'\^etre for the strongly stable and the
    stable partial orders: Borel ideals and stable ideals}

  \begin{subsection}{Borel ideals and stable ideals}
  \begin{definition}\label{def:Borel}
    Let \(n,d\) be positive integers. A subset \(U \subset \M_d\) is
    called \emph{strongly stable} or \emph{Borel} iff
    \begin{displaymath}
      m \in U,  \, \divides{x_j}{m}, \,  1 \le i \le j
      \le n \quad \implies \quad m \frac{x_i}{x_j} \in U.
    \end{displaymath}
    A monoid ideal \(I \subset \M\) is called a strongly stable
    monoid ideal or a
    Borel monoid ideal iff \(I \cap \M_v\) is Borel for all
    positive integers \(v\).

    Borel subsets of \(M_d^n\) and Borel ideals in \(\M^n\) are
    defined 
    analogously. 

    If \(K\) is a field, then a monomial ideal \(J \subset K\M \simeq
    K[x_1,x_2,x_3,\dots]\) is called Borel (or strongly stable) if \(J
    \cap \M\) is a Borel monoid ideal, and similarly for monomial
    ideals in \(K\M^n \simeq K[x_1,\dots,x_n]\).
  \end{definition}

  A reason to study Borel ideals is 
  following theorem (see \cite{Galligo:Prep, Galligo:DivStab,
  Green:gin, Ebud:View, BS:Refine}). Recall that a \emph{term order} on \(\M^n\)
  is a strongly multiplicative total order \(\le\) such that  \(1\) is the
  smallest 
  element. Given a term order \(\le\) and a polynomial \(f \in \Kxn\),
  the \emph{initial monomial} \(\init_\le(f)\) is the largest (
  \wrt \(\le\)) monomial occurring in \(f\). If \(J \subset
  \Kxn\) is an ideal, then the \emph{initial ideal} \(\init_\le(J)\)
  is the ideal generated by 
  \(\setsuchas{\init_\le(f)}{f \in J}\). 
  Finally, the general linear group \(\mathrm{GL}_n\) acts on
  the \(n\)-dimensional \(K\)-vector space \(V\) spanned by the
  variables in \(\Kxn\), and since \(\Kxn\) is the symmetric algebra
  on \(V\), this action can be extended to \(\Kxn\). Explicitly, if
  \(\mathrm{GL}_n \ni g = (g_{ij})\) then \(g\) acts on the monomial
  \(x_1^{\alpha_1} \cdots x_n^{\alpha_n}\) by 
  \begin{displaymath}
    g\left(\prod_{i=1}^n x_i^{\alpha_i}\right) = \prod_{i=1}^n
    g(x_i)^{\alpha_i} 
    = \prod_{i=1}^n {\left(\sum_{j=1}^n g_{ij}x_j \right)}^{\alpha_i},    
  \end{displaymath}
  and this action is then extended \(K\)-linearly.

  \begin{theorem}[Galligo, Bayer-Stillman]\label{thm:galligo}
    Let \(\ge\) be a term order with \(x_1 \ge \cdots \ge x_n\), let
    \(K\) be a  
    field of characteristic 0, let \(n\) be a positive integer, and
    let \(J \subset \Kxn\) be a homogeneous ideal. Then there is a
    monomial ideal 
    \(\mathrm{gin}_\ge(J)\), the \emph{generic initial ideal} of
    \(J\), and a Zariski-open subset \(U \subset \mathrm{GL}_n\), such
    that \(\mathrm{in}(g(J)) = \mathrm{gin}_\ge(J)\) for 
    all \(g \in U\).
    Furthermore \(\mathrm{gin}_\ge(J) \cap 
    \M^n\) is   a Borel monoid ideal.
  \end{theorem}

  \begin{definition}\label{def:Stabideal}
    Let \(n,d\) be positive integers. A subset \(U \subset \M_d\) is
    called \emph{stable} iff
    \begin{displaymath}
      m \in U,  \,   i \le \maxsupp(m) \quad \implies \quad m
      \frac{x_i}{x_{\maxsupp(m)}} \in U. 
    \end{displaymath}
    A monoid ideal \(I \subset \M\) is called a stable
    monoid ideal iff \(I \cap \M_v\) is Borel for all
    positive integers \(v\).

    Stable subsets of \(M_d^n\) and stable ideals in \(\M^n\) are
    defined 
    analogously. 

    If \(K\) is a field, then a monomial ideal \(J \subset K\M \simeq
    K[x_1,x_2,x_3,\dots]\) is called  stable if \(J
    \cap \M\) is a stable monoid ideal, and similarly for monomial
    ideals in \(K\M^n \simeq K[x_1,\dots,x_n]\).
  \end{definition}

  Stable ideals have minimal free resolutions of a particularly nice
  form \cite{Eliahou:MinRes}, and are therefore interesting.

  The following theorem motivates the study of the
  stable and strongly stable partial orders:
  \begin{thm}\label{thm:filters}
    Let \(d,n\) be positive integers.
    \begin{enumerate}[(i)]
    \item \label{enum:f1}
      A subset \(I \subset \M\) is a monoid ideal iff it is a
      filter \wrt the partial order \(D\).
    \item \label{enum:f2}
      A subset \(I \subset \M^n\) is a monoid ideal iff it is a
      filter \wrt the partial order \(D\).
    \item \label{enum:f3}
      A subset \(U \subset \M_d\) is  Borel  iff it is a filter
      \wrt the partial order \(A_{\cdot,d}\).
    \item \label{enum:f4}
      A subset \(U \subset \M_d^n\) is  Borel  iff it is a filter
      \wrt the partial order \(A_{n,d}\).
    \item \label{enum:f5}
      A subset \(I \subset \M\) is a Borel monoid ideal iff it is
      a filter 
      \wrt the partial order \(A_{\cdot,\cdot}\).
    \item \label{enum:f6}
      A subset \(I \subset \M^n\) is a Borel monoid ideal iff it
      is a filter 
      \wrt the partial order \(A_{n,\cdot}\).
    \item \label{enum:f7}
      A subset \(I \subset \M\) is a stable monoid ideal iff it is
      a filter 
      \wrt the partial order \(B_{\cdot,\cdot}\).
    \item \label{enum:f8}
      A subset \(I \subset \M^n\) is a stable monoid ideal iff it
      is a filter 
      \wrt the partial order \(B_{n,\cdot}\).
    \end{enumerate}
  \end{thm}
  \begin{proof}
    \prettyref{enum:f1} and \prettyref{enum:f2} are well-known, and 
    \prettyref{enum:f3} and \prettyref{enum:f4}
    is immediate from the definitions. \prettyref{enum:f5},
    \prettyref{enum:f6},  \prettyref{enum:f7}, and \prettyref{enum:f8}
     are similar; we prove \prettyref{enum:f5}.

    If \(I\) is a filter \wrt \(A_{\cdot,\cdot}\), then
    since \(D \subset A_{\cdot,\cdot}\),
    it is a 
    filter \wrt \(D\), hence \(I\) is a monoid
    ideal. For any \(d \in \Nat^+\) we have that \(I \cap \M_d\) is a
    filter \wrt \(A_{\cdot,d}\), since
    \(\restRel{A_{\cdot,\cdot}}{\M_d} =  A_{\cdot,d}\). This
    shows that \(I_d\) is Borel. Hence, \(I\) is a Borel monoid ideal.

    Conversely, if \(I\) is a Borel monoid ideal, we want to show that
    \(I\) is a filter \wrt  \(A_{\cdot,\cdot}\). Since it
    is a monoid ideal, it is a filter \wrt \(D\). By
    definition, for each \(d\) we have that \(I_d\)  is a filter \wrt
     \(A_{\cdot,d}\). Since \(A_{\cdot,\cdot}\) is the
    smallest partial order which contains \(D\) and all
    \(A_{\cdot,d}\), it follows that \(I\) is a filter \wrt
    \(A_{\cdot,\cdot}\). 
  \end{proof}
  \end{subsection}

\end{section}

\begin{section}{Properties of the strongly stable partial order}
    \begin{subsection}{Properties of infinite strongly stable partial orders}
  We denote the dual of any partial order \(P\) by
  \(\dual{P}\). 
  \begin{definition}
    For \(n,d \in \Nat\) we define 
    \begin{align}
    C_{n,d} &= \dual{A_{n,d}}\\
    C_{\cdot,d} &= \bigcup_{n} C_{n,d} \subset \M_d \times \M_d\\
    C_{n,\cdot} &= \rtr(D \cup \bigcup_{d} C_{n,d}) \subset \M^n \times \M^n
    \\
    C_{\cdot,\cdot} & := \rtr(D \cup \bigcup_{n,d} C_{n,d})
    \end{align}
  \end{definition}
  Clearly, \(C_{n,\cdot} \simeq A_{n,\cdot}\) for any \(n\), because
  we can simply rename the variables according to the bijection \(x_i
  \leftrightarrow x_{n+1-i}\). However, \(A_{\cdot,\cdot} \not
  \simeq C_{\cdot,\cdot}\):  \(C_{\cdot,\cdot}\) has a smallest element,
  \(1\), and an element covering it, \(x_1\).
  \(A_{\cdot,\cdot}\), on the other hand, has a smallest element, \(1\), but
  no element covering it.

  We denote by  \(\mathcal{Y}\) the Young lattice of decreasing,
  eventually zero sequences of non-negative integers, ordered by
  component-wise inclusion (see \cite{Aigner:Combinatorial} for a more
  thorough treatment).

  \begin{theorem}\label{thm:CYoung}
    \(C_{\cdot,\cdot} \simeq \mathcal{Y}\).
  \end{theorem}
  \begin{proof}
    To the monomial \(x_1^{\alpha_1} \cdots x_n^{\alpha_n}\) we
    associate the Ferrers diagram consisting of \(\alpha_n\) rows of
    length \(n\), \(\alpha_{n-1}\) rows of length \(n-1\), and so
    on. This is clearly a bijection between \(\M\) and
    \(\mathcal{Y}\). We now show that it is isotone.

    Consider all relations of the following two types: 
    \begin{enumerate}
    \item \(x_i m \mapsto x_{i+1} m\),
    \item \(m \mapsto x_i m\).
    \end{enumerate}
    These two relations generate the partial order.
    One sees that type 1 corresponds to enlarging the topmost row of
    the rows with \(i\) element with 1 element, and that type 2
    corresponds to inserting a row with \(i\) elements. Therefore, the
    map is isotone. Furthermore, since these two types of operations
    on Ferrers diagrams generate the Young lattice, the inverse is
    isotone as well.
  \end{proof}

  \begin{corr}
    \(C_{\cdot,\cdot}\) is a distributive lattice. For any \(n,d
    \in \Nat\), the posets \(A_{n,\cdot} \simeq
    \dual{C_{n,\cdot}}\) and \(A_{n,d} \simeq
    \dual{C_{n,d}}\)  
    are distributive lattices, isomorphic to the set of Ferrers
    diagrams with at most \(n\) columns, or with at most \(n\) columns
    and exactly \(d\) rows, respectively.
  \end{corr}

  \begin{prop}\label{prop:remove}
    For any positive integers \(v,w\), the subset of all Ferrers
    diagrams 
    with exactly \(v\) rows is a poset isomorphic to the set of all
    Ferrers diagrams with at most \(v\) rows. Furthermore, the subset
    of  all Ferrers diagrams with at most \(w\) columns and exactly
    \(v\) rows is isomorphic to the set of all Ferrers diagrams which
    fit inside a \(v \times (w-1)\) box.
  \end{prop}
  \begin{proof}
    If a Ferrers diagram has exactly \(v\) rows, then its first column
    has length \(v\). Removing this column, one gets a Ferrers
    diagrams with one less column, and with at most \(v\) rows. This
    establishes the desired bijections, which are isotone
    with isotone inverses.
  \end{proof}

  The situation for \(A_{\cdot,\cdot}\) is different, but
  similar. 
  \begin{theorem}\label{thm:Adist}
    The poset \((\M,A_{\cdot, \cdot})\) is a distributive
    lattice, isomorphic to the set of all weakly increasing,
    eventually constant functions \(\Nat^+ \to \Nat\), ordered by
    component-wise inclusion.
  \end{theorem}
  \begin{proof}[Sketch of proof]
    Use the map 
    \begin{align*}
      \Integ: \,M & \to \Nat^{\Nat^+}\\
      x_1^{\alpha_1} \cdots x_n^{\alpha_n} & \mapsto
      (\alpha_1, \alpha_1 + \alpha_2, \dots, \alpha_1 +\dots+\alpha_n,
      \alpha_1 +\dots +\alpha_n, \dots)
    \end{align*}
    It is not hard to see that \(\Xi\) is injective, and that its image is all
    weakly increasing sequences which are eventually
    constant. Furthermore, one can convince oneself that it is
    isotone: the operation \(m \to x_i m\) maps to the
    operation of adding the sequence \(B[i]\) which is 1 from \(i\)
    and onward 
    and zero before that, and the operation \(m \to (x_i/x_{i+1})m\)
    maps to the operation of inserting a block of height and width 1,
    at position \(i+1\). Formally, the sequence
    \((\dots,b_i,b_{i+1}, b_{i+2},\dots)\), which is required to have
    a ``jump'' between \(i\) and \(i+1\), that is, \(b_i < b_{i+1}\),
    is replaced with the 
    sequence \((\dots, 1 + b_i, b_{i+1}, b_{i+2},\dots)\).
    
    The hard part is showing that the inverse \(\Deriv\) is isotone.
    For this, one need to show that the two operations of adding
    sequences of the form \(B[i]\), and inserting a single block at a
    ``jump'', generates the order relation for weakly increasing,
    eventually constant sequences. One can prove this by induction
    over the eventually constant value, and over the point from which
    it becomes constant.
  \end{proof}

  \begin{corr}
    For any \(n,d \in \Nat\), the posets \(A_{n,\cdot}\) and
    \(A_{\cdot,d}\) are
    distributive sublattice of the distributive lattice
    \(A_{\cdot,\cdot}\). They are isomorphic to the following
    two subsets of the set of all weakly increasing,
    eventually constant functions \(\Nat^+ \to \Nat\), ordered by point-wise
    comparison:  
    \begin{itemize}
    \item The set \(\Integ(\M^n)\) of such functions \(f\) with
      \(f(n) = \lim_{t \to +\infty} f(t) \),
    \item The set \(\Integ(\M_d)\) of such functions \(f\) with
      \(d= \lim_{t \to +\infty} f(t)\).
    \end{itemize}
  \end{corr}
  
  Another result that follows immediately from the above description is the
  following: 
  \begin{corr}\label{corr:antichainprod}
    If \(x,y \in \M\) and \(\set{x,y}\) is an antichain \wrt
    \(A_{\cdot,\cdot}\) then 
    so is \(\set{x^n,y^n}\), for all \(n \in \Nat^+\).
  \end{corr}
  \begin{proof}
   By the previous theorem this is translated into the following
    assertion: if \(f,g \in \Integ(\M)\) are incomparable, then so is
    \(nf\) and \(ng\). Now, \(f,g\)
    are incomparable iff there exists \(a,b \in \Nat^+\) with \(f(a) >
    g(a)\), \(f(b) < g(b)\), and this implies that 
    \(nf(a) >  ng(a)\), \(nf(b) < ng(b)\), showing that \(nf\)
    and \(ng\) are incomparable.
  \end{proof}
  

  \begin{prop}\label{prop:anti}
    For any positive integer \(v\), 
    \((\M^v,\,A_{v,\cdot})\)
    and 
    \((\M_v,\, A_{\cdot,v})\)
    are dual posets.
  \end{prop}
  \begin{proof}
    It will suffice to give an antitone bijection \(\tau\) with
    antitone 
    inverse \(\tau^{-1}\) between \(\M^v\) and
    \(\Integ(\M_v)\). This map is defined as follows. 
    Take \[\hat{a} =(a_1,\dots,a_v,0,0,\dots) \in \M^v.\] 
    Put \(\ell =   \tdeg{\hat{a}} = \sum_{i=1}^v a_i\). Define
    \[\tau(\hat{a}) = (b_1,\dots,b_\ell,b_\ell,\dots) \in \Integ(\M_v),\]
    with
    \begin{displaymath}
      b_i = 
      \begin{cases}
        0 &  i \le  a_1\\
        1 & a_1 < i \le a_1 + a_2 \\
        \vdots  & \vdots\\
        v-1 & a_1 + \cdots + a_{v-1} < i \le a_1 + \cdots + a_v \\
        v & i > \ell
      \end{cases}
    \end{displaymath}
    That is, \(\tau(\hat{a})\) has \(a_1\) zeroes, \(a_2\) ones, and
    so on, and takes on the constant value \(v\) at position \(\ell\) and
    onward.
    It is easily seen that \(\tau\) is bijective, and it is fact also
    antitone with antitone inverse.
  \end{proof}
  
  For \(v=2\), the (beginning of) the Hasse diagrams for these two
  infinite posets are depicted in \prettyref{fig:2varHasse}.
  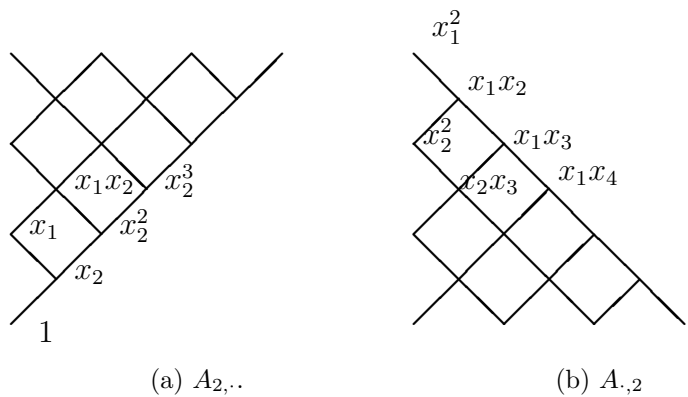
\begin{figure}[tbh]
    \begin{center}
      \setlength{\unitlength}{1.2cm}
      \subfigure[\(A_{2,\cdot}\).]{
\begin{picture}(4,4) \thicklines
\put(0,0){\line(1,1){3}}
\put(0,1){\line(1,1){2}}
\put(0,1){\line(1,-1){0.5}}
\put(0,2){\line(1,1){1}}
\put(0,2){\line(1,-1){1}}
\put(0,3){\line(1,-1){1.5}}
\put(1,3){\line(1,-1){1}}
\put(2,3){\line(1,-1){0.5}}
\put(0.3,-0.2){\(1\)}
\put(0.7,0.5){\(x_2\)}
\put(0.2,1){\(x_1\)}
\put(1.2,1){\(x_2^2\)}
\put(0.7,1.5){\(x_1x_2\)}
\put(1.7,1.5){\(x_2^3\)}
\end{picture}
}      
\subfigure[\(A_{\cdot,2}\)]{
\begin{picture}(4,4) \thicklines
\put(0,3){\line(1,-1){3}}
\put(0,2){\line(1,-1){2}}
\put(0,2){\line(1,1){0.5}}
\put(0,1){\line(1,-1){1}}
\put(0,1){\line(1,1){1}}
\put(0,0){\line(1,1){1.5}}
\put(1,0){\line(1,1){1}}
\put(2,0){\line(1,1){0.5}}
\put(0.2,3.2){\(x_1^2\)}
\put(0.6,2.6){\(x_1x_2\)}
\put(0.1,2.0){\(x_2^2\)}
\put(1.1,2.0){\(x_1x_3\)}
\put(0.5,1.5){\(x_2x_3\)}
\put(1.6,1.6){\(x_1x_4\)}
\end{picture}
        }
    \end{center}
    \caption{The Hasse diagrams for the strongly stable order on 2
      variables, or on monomials of degree 2.}
      \label{fig:2varHasse}
    \end{figure}
  \end{subsection}

  \begin{subsection}{Properties of finite strongly stable partial
    orders} 
    We fix positive integers \(n,d\), and study \(A_{n,d}\) and
    \(A_{n,d}^\circ\). 

    We note that \(A_{1,d}\) is a singleton for all \(d\), and 
    that \(A_{2,d}\) is a chain for all \(d\). 
    The relation \(A_{3,d}^\circ\) and the Hasse diagrams of 
    \(A_{3,d}\) looks as 
    \prettyref{fig:HasseSS3}.


  We recall that
  \(C_{n,d} \simeq A_{n,d}\), so that
    \(A_{n,d}\) is self-dual and isomorphic with the set of
    all Ferrers diagrams which fits inside a \(d \times (n-1)\) box, ordered
    by inclusion.

    \begin{definition}\label{def:qbinom}
      The \(q\)-binomial, or Gaussian polynomials are defined as
      \begin{multline}
        {\left[
            \begin{array}[c]{c}
              a+b \\
              a
            \end{array}
          \right]}_q =
        \frac{(1-q^{a+b})(1-a^{a+b-1}) \cdots
          (1-q^{a+1})}{(1-q^b)(1-q^{b-1})\cdots (1-q)} = \\ =
        c_0^{(a,b)} + c_1^{(a,b)}q + \cdots + c_N^{(a,b)}q^N \in \Z[q].
      \end{multline}
      The coefficients \(c_{i}^{(a,b)}\) are called the \(q\)-binomial 
      coefficients.
    \end{definition}
    
    We list some well-known properties of the poset \(F_{a,b} \simeq
    (\M_{b}^{a+1}, \rtr(A_{a+1,b}))\) of all Ferrers
    diagrams that fit inside an \(a \times b\) rectangle.
    For the definition of the height, width and dimension of a
    finite partially ordered set, see \cite{Trotter:Poset}.

      \begin{itemize}
      \item The poset \(F_{a,b}\) is a ranked distributive lattice,
        with rank function \(f(i) = c_{i}^{(a,b)}\).
        \item The ranked poset \(F_{a,b}\) is rank unimodal, rank
          symmetric,  and Sperner.
        \item \(F_{a,b}\) has dimension \(a\), height \(ab\), and
          width \(c_v^{(a,b)}\), where \(v={\lfloor ab/2 \rfloor}\).
      \end{itemize}
    
  \begin{corr}
    For positive integers \(n,d \ge 2\), the poset \(A_{n,d}\)
    has dimension \(n-1\), height \((n-1)d\), and width
    \(c_w^{(n-1,d)}\), where \(w=\lfloor (n-1)d/2 \rfloor\).
  \end{corr}

  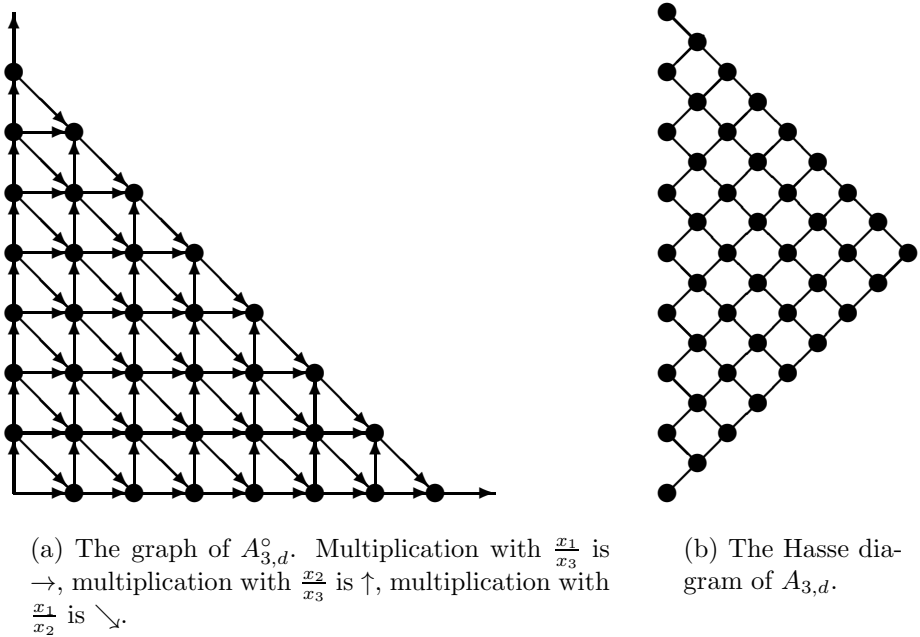
\begin{figure}[tbh]
    \begin{center}
      \setlength{\unitlength}{0.8cm}
      \subfigure[The graph of \(A_{3,d}^\circ\).
        Multiplication
        with \(\frac{x_1}{x_3}\) is 
        \(\rightarrow\), multiplication with \(\frac{x_2}{x_3}\) is
        \(\uparrow\), 
        multiplication with \(\frac{x_1}{x_2}\) is \(\searrow\).
        ]{
            \begin{picture}(10,9)    \thicklines
    \put(0,0){\vector(1,0){8}}
    \put(0,0){\vector(0,1){8}}
    \put(1,0){\circle*{0.3}}
    \put(0,1){\circle*{0.3}}
    \put(2,0){\circle*{0.3}}
    \put(1,1){\circle*{0.3}}
    \put(0,2){\circle*{0.3}}
    \put(3,0){\circle*{0.3}}
    \put(2,1){\circle*{0.3}}
    \put(1,2){\circle*{0.3}}
    \put(0,3){\circle*{0.3}}
    \put(4,0){\circle*{0.3}}
    \put(3,1){\circle*{0.3}}
    \put(2,2){\circle*{0.3}}
    \put(1,3){\circle*{0.3}}
    \put(0,4){\circle*{0.3}}
    \put(5,0){\circle*{0.3}}
    \put(4,1){\circle*{0.3}}
    \put(3,2){\circle*{0.3}}
    \put(2,3){\circle*{0.3}}
    \put(1,4){\circle*{0.3}}
    \put(0,5){\circle*{0.3}}
    \put(6,0){\circle*{0.3}}
    \put(5,1){\circle*{0.3}}
    \put(4,2){\circle*{0.3}}
    \put(3,3){\circle*{0.3}}
    \put(2,4){\circle*{0.3}}
    \put(1,5){\circle*{0.3}}
    \put(0,6){\circle*{0.3}}
    \put(7,0){\circle*{0.3}}
    \put(6,1){\circle*{0.3}}
    \put(5,2){\circle*{0.3}}
    \put(4,3){\circle*{0.3}}
    \put(3,4){\circle*{0.3}}
    \put(2,5){\circle*{0.3}}
    \put(1,6){\circle*{0.3}}
    \put(0,7){\circle*{0.3}}
    \put(0,0){\vector(1,0){0.9}}
    \put(0,0){\vector(0,1){0.9}}
    \put(1,0){\vector(1,0){0.9}}
    \put(1,0){\vector(0,1){0.9}}
    \put(0,1){\vector(1,0){0.9}}
    \put(0,1){\vector(0,1){0.9}}
    \put(2,0){\vector(1,0){0.9}}
    \put(2,0){\vector(0,1){0.9}}    
    \put(1,1){\vector(1,0){0.9}}
    \put(1,1){\vector(0,1){0.9}}
    \put(0,2){\vector(1,0){0.9}}
    \put(0,2){\vector(0,1){0.9}}
    \put(3,0){\vector(1,0){0.9}}
    \put(3,0){\vector(0,1){0.9}}    
    \put(2,1){\vector(1,0){0.9}}
    \put(2,1){\vector(0,1){0.9}}    
    \put(1,2){\vector(1,0){0.9}}
    \put(1,2){\vector(0,1){0.9}}    
    \put(0,3){\vector(1,0){0.9}}
    \put(0,3){\vector(0,1){0.9}}    
    \put(4,0){\vector(1,0){0.9}}
    \put(4,0){\vector(0,1){0.9}}    
    \put(3,1){\vector(1,0){0.9}}
    \put(3,1){\vector(0,1){0.9}}    
    \put(2,2){\vector(1,0){0.9}}
    \put(2,2){\vector(0,1){0.9}}    
    \put(1,3){\vector(1,0){0.9}}
    \put(1,3){\vector(0,1){0.9}}    
    \put(0,4){\vector(1,0){0.9}}
    \put(0,4){\vector(0,1){0.9}}    
    \put(5,0){\vector(1,0){0.9}}
    \put(5,0){\vector(0,1){0.9}}    
    \put(4,1){\vector(1,0){0.9}}
    \put(4,1){\vector(0,1){0.9}}    
    \put(3,2){\vector(1,0){0.9}}
    \put(3,2){\vector(0,1){0.9}}    
    \put(2,3){\vector(1,0){0.9}}
    \put(2,3){\vector(0,1){0.9}}    
    \put(1,4){\vector(1,0){0.9}}
    \put(1,4){\vector(0,1){0.9}}    
    \put(0,5){\vector(1,0){0.9}}
    \put(0,5){\vector(0,1){0.9}}    
    \put(6,0){\vector(1,0){0.9}}
    \put(6,0){\vector(0,1){0.9}}    
    \put(5,1){\vector(1,0){0.9}}
    \put(5,1){\vector(0,1){0.9}}    
    \put(4,2){\vector(1,0){0.9}}
    \put(4,2){\vector(0,1){0.9}}    
    \put(3,3){\vector(1,0){0.9}}
    \put(3,3){\vector(0,1){0.9}}    
    \put(2,4){\vector(1,0){0.9}}
    \put(2,4){\vector(0,1){0.9}}    
    \put(1,5){\vector(1,0){0.9}}
    \put(1,5){\vector(0,1){0.9}}    
    \put(0,6){\vector(1,0){0.9}}
    \put(0,6){\vector(0,1){0.9}}    
    \put(0,7){\vector(1,-1){0.9}}
    \put(1,6){\vector(1,-1){0.9}}
    \put(2,5){\vector(1,-1){0.9}}
    \put(3,4){\vector(1,-1){0.9}}
    \put(4,3){\vector(1,-1){0.9}}
    \put(5,2){\vector(1,-1){0.9}}
    \put(6,1){\vector(1,-1){0.9}}
    \put(0,6){\vector(1,-1){0.9}}
    \put(1,5){\vector(1,-1){0.9}}
    \put(2,4){\vector(1,-1){0.9}}
    \put(3,3){\vector(1,-1){0.9}}
    \put(4,2){\vector(1,-1){0.9}}
    \put(5,1){\vector(1,-1){0.9}}
    \put(0,5){\vector(1,-1){0.9}}
    \put(1,4){\vector(1,-1){0.9}}
    \put(2,3){\vector(1,-1){0.9}}
    \put(3,2){\vector(1,-1){0.9}}
    \put(4,1){\vector(1,-1){0.9}}
    \put(0,4){\vector(1,-1){0.9}}
    \put(1,3){\vector(1,-1){0.9}}
    \put(2,2){\vector(1,-1){0.9}}
    \put(3,1){\vector(1,-1){0.9}}
    \put(0,3){\vector(1,-1){0.9}}
    \put(1,2){\vector(1,-1){0.9}}
    \put(2,1){\vector(1,-1){0.9}}
    \put(0,2){\vector(1,-1){0.9}}
    \put(1,1){\vector(1,-1){0.9}}
    \put(0,1){\vector(1,-1){0.9}}
  \end{picture}
          }      
        \setlength{\unitlength}{0.4cm}
        \subfigure[The Hasse diagram of \(A_{3,d}\).
        ]{
          \begin{picture}(8,15) \thicklines
\put(1,1){\line(-1,1){1.2}}
\put(2,2){\line(-1,1){2}}
\put(3,3){\line(-1,1){3}}
\put(4,4){\line(-1,1){4}}
\put(5,5){\line(-1,1){5}}
\put(6,6){\line(-1,1){6}}
\put(7,7){\line(-1,1){7}}
\put(8,8){\line(-1,1){8}}
\put(0,0){\line(1,1){8}}
\put(0,2){\line(1,1){7}}
\put(0,4){\line(1,1){6}}
\put(0,6){\line(1,1){5}}
\put(0,8){\line(1,1){4}}
\put(0,10){\line(1,1){3}}
\put(0,12){\line(1,1){2}}
\put(0,14){\line(1,1){1.2}}
\put(0,0){\circle*{0.6}}
\put(0,2){\circle*{0.6}}
\put(0,4){\circle*{0.6}}
\put(0,6){\circle*{0.6}}
\put(0,8){\circle*{0.6}}
\put(0,10){\circle*{0.6}}
\put(0,12){\circle*{0.6}}
\put(0,14){\circle*{0.6}}
\put(0,16){\circle*{0.6}}
\put(1,1){\circle*{0.6}}
\put(1,3){\circle*{0.6}}
\put(1,5){\circle*{0.6}}
\put(1,7){\circle*{0.6}}
\put(1,9){\circle*{0.6}}
\put(1,11){\circle*{0.6}}
\put(1,13){\circle*{0.6}}
\put(1,15){\circle*{0.6}}
\put(2,2){\circle*{0.6}}
\put(2,4){\circle*{0.6}}
\put(2,6){\circle*{0.6}}
\put(2,8){\circle*{0.6}}
\put(2,10){\circle*{0.6}}
\put(2,12){\circle*{0.6}}
\put(2,14){\circle*{0.6}}
\put(3,3){\circle*{0.6}}
\put(3,5){\circle*{0.6}}
\put(3,7){\circle*{0.6}}
\put(3,9){\circle*{0.6}}
\put(3,11){\circle*{0.6}}
\put(3,13){\circle*{0.6}}
\put(4,4){\circle*{0.6}}
\put(4,6){\circle*{0.6}}
\put(4,8){\circle*{0.6}}
\put(4,10){\circle*{0.6}}
\put(4,12){\circle*{0.6}}
\put(5,5){\circle*{0.6}}
\put(5,7){\circle*{0.6}}
\put(5,9){\circle*{0.6}}
\put(5,11){\circle*{0.6}}
\put(6,6){\circle*{0.6}}
\put(6,8){\circle*{0.6}}
\put(6,10){\circle*{0.6}}
\put(7,7){\circle*{0.6}}
\put(7,9){\circle*{0.6}}
\put(8,8){\circle*{0.6}}
\end{picture}
          }
    \end{center}
    \caption{The strongly stable relation and partial order for
      \(n=3\).} 
      \label{fig:HasseSS3}
    \end{figure}

  \end{subsection}

  \begin{subsection}{Relation to term orders}
    We mean by a \emph{term order} a strongly multiplicative total
    order on \(\M\) or on \(\M^n\), with \(1 < m\) for all \(m \neq 1\). 
  \begin{definition}\label{def:totorders}
    Let \(n\) be a positive integer.     
    Denote by \(\mathfrak{T}\) the set of all term orders
    orders \(R\) on \(\M\) such that \(x_1 \,R \, x_2\,R \,x_3 \, R
    \cdots\).
    Similarly, denote by \(\mathfrak{U}\) the set of all 
    term orders  \(S\) on \(\M^n\) such that \(x_1 \,S \, x_2\,S \,x_3
    \, S \cdots S x_n\).
  \end{definition}

  \begin{theorem}\label{thm:termorderIntersect}
    Let \(n\) be a positive integer.
    Then we have that 
    \begin{equation}\label{eqn:inters}
      \begin{split}
      A_{\cdot,\cdot} &= \bigcap_{R \in \mathfrak{T}} R \\
      A_{n,\cdot} &= \bigcap_{S\in \mathfrak{U}} S \\
      \end{split}
    \end{equation}
  \end{theorem}

  \begin{proof}
    To start, note that \(\bigcap_{R \in \mathfrak{T}} R\) and
    \(\bigcap_{S\in \mathfrak{U}} S\) are posets.
    Take
    \((m,m') \in A_{\cdot,\cdot}^\circ\). If \(\divides{m'}{m}\) then 
    \((m,m') \in R\)
    for every \(R \in \mathfrak{T}\), since such an \(R\) is
    multiplicative. If on the other hand \((m,m') \in 
    A_{\cdot,\cdot}^\circ\) but \(\dividesnot{m'}{m}\) then we may assume
    that \(m = \frac{x_i}{x_j} m'\),    with \(i <    j\). Any \(R \in
    \mathfrak{T}\) may be 
    extended to a     multiplicative total order on the difference
    group of \(\M\), and     it 
    is clear that for this extension we have that \((\frac{x_i}{x_j},1) \in
    R\), 
    hence \((m' \frac{x_i}{x_j},m') \in
    R\), since \(R\) is multiplicative. Since \(D \subset R\) and \(A_{n,d}
    \subset R\) for all 
    \(n,d\), and since \(R\) is a strongly multiplicative total order,
    we have that \(A_{\cdot,\cdot} \subset R\). Since \(R\) was
    arbitrary, we have
    proved that 
    \(A_{\cdot,\cdot} \subset \left(\bigcap_{R \in   \mathfrak{T}} R\right) .\)
  Hence, every multiplicative total order on \(\M\) extends
  \(\rtr(A_{\cdot,\cdot})\). 
  
  Suppose now that \((m,m') \not \in A_{\cdot,\cdot}\). We must show
  that \((m,m') \not \in \bigcap_{R \in \mathfrak{T}} R\). To do this,  we
  show that there exist
  some term order \(R\) such that \((m',m) \in R\). If \((m',m) \in
  A_{\cdot,\cdot}\) then the argument above shows that \((m',m) \in
  R\) for all \(R \in \mathfrak{T}\). We address the remaining case,
  where \(m,m'\) are incomparable. Then from \prettyref{corr:antichainprod} we 
  have that \(m^n\) and \({m'}^n\) are incomparable for all \(n \in \Nat^+\),
  thus from \prettyref{thm:gilmerplus} we get 
  that \(A_{\cdot,\cdot}\) can be extended to a multiplicative total order
  \(R\) such that \((m',m) \in R\).
  \end{proof}

  So every  term order \(>\) satisfying \(x_1 > x_2 >\cdots\)
  refines the strongly stable partial order, and \((m, m') \in
  \rtr(A_{\cdot,\cdot})\) iff \(m \ge m'\) for all such admissible orders.

  \begin{corr}\label{corr:filters}
    Let \(R \in \mathfrak{T}\)   and let \(U \subset \M\) be a filter
    \wrt \(R\). Then  \(U\) is a filter \wrt
    \(A_{\cdot,\cdot}\) (and is thus a Borel monoid ideal in
    \(\M\)). 

    If \(n\) is a positive integer, if \(S
    \in \mathfrak{U}\), and if \(V \subset \M^n\) is a filter \wrt
     \(S\), then \(V\) is a filter \wrt
    \(A_{n,\cdot}\) (and is thus a Borel monoid ideal in
    \(\M^n\)). 
  \end{corr}

  We call a term-order \(\ge\) \emph{degree-compatible} if it refines
  the   partial order given by total degree: in other words, if
  \(\tdeg{m} > \tdeg{m'} \implies m > m'\). 
  
  \begin{definition}\label{def:degtotorders}
    Let \(n\) be a positive integer.     
    Denote by \(\mathfrak{dT}\) the set of all degree-compatible
    term orders \(R\) on \(\M\) such that 
    \(x_1 \,R \, x_2\,R \,x_3 \, R \cdots\). 
    Similarly, denote by \(\mathfrak{dU}\) the set of all
    degree-compatible term orders \(S\) on \(\M^n\) such that \(x_1
    \,S \, x_2\,S \,x_3 \, S \cdots S x_n\). 
  \end{definition}
  
  \begin{theorem}\label{thm:degcomptermorderIntersect}
    Let \(n\) be a positive integer.
    Then we have that 
    \begin{equation}\label{eqn:dinters}
      \begin{split}
      \bigoplus_{d=0}^\infty A_{\cdot,d} &= \bigcap_{R \in
      \mathfrak{dT}} R \\ 
      \bigoplus_{d=0}^\infty A_{n,d} &= \bigcap_{S\in
      \mathfrak{dU}} S \\ 
      \end{split}
    \end{equation}
    Here, \(\bigoplus\) denotes the ordinal sum of posets (see
    \cite[VIII,\S 10]{Birkhoff:LT}).
  \end{theorem}
  \begin{proof}
    Similar to \prettyref{thm:termorderIntersect}, noting that there
    can be no antichain between monomials of different total degree.
  \end{proof}


  \begin{figure}[htbp]
  \begin{center}
    \leavevmode
    \setlength{\unitlength}{0.3cm}
        \begin{picture}(8,21) \thicklines
\put(0,0){\circle*{0.2}}
\put(0,0){\line(0,1){3}}

\put(0,3){\circle*{0.2}}
\put(2,4){\circle*{0.2}}
\put(0,5){\circle*{0.2}}
\put(0,3){\line(2,1){2}}
\put(2,4){\line(-2,1){2}}
\put(0,5){\line(0,1){3}}

\put(0,8){\circle*{0.2}}
\put(0,10){\circle*{0.2}}
\put(0,12){\circle*{0.2}}
\put(2,9){\circle*{0.2}}
\put(2,11){\circle*{0.2}}
\put(4,10){\circle*{0.2}}
\put(0,8){\line(2,1){4}}
\put(0,10){\line(2,1){2}}
\put(2,9){\line(-2,1){2}}
\put(4,10){\line(-2,1){4}}
\put(0,12){\line(0,1){3}}

\put(0,15){\circle*{0.2}}
\put(0,17){\circle*{0.2}}
\put(0,19){\circle*{0.2}}
\put(0,21){\circle*{0.2}}
\put(2,16){\circle*{0.2}}
\put(2,18){\circle*{0.2}}
\put(2,20){\circle*{0.2}}
\put(4,17){\circle*{0.2}}
\put(4,19){\circle*{0.2}}
\put(6,18){\circle*{0.2}}
\put(0,15){\line(2,1){6}}
\put(2,16){\line(-2,1){2}}
\put(4,17){\line(-2,1){4}}
\put(6,18){\line(-2,-1){6}}
\put(6,18){\line(-2,1){6}}
\put(4,19){\line(-2,-1){4}}
\put(2,20){\line(-2,-1){2}}

\end{picture}
  \caption{Hasse diagram of \(\bigoplus_{d=1}^\infty A_{3,d}\).}
    \label{fig:ordsum}
  \end{center}
\end{figure}
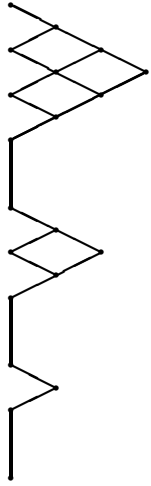

  The Hasse diagram for
  \(\bigoplus_{d=1}^\infty A_{2,d}\) is the ordinal sum of
  chains, hence a chain. The Hasse diagram for
  \(\bigoplus_{d=1}^\infty A_{3,d}\) is more interesting: it
  looks like figure~\ref{fig:ordsum}.
  \end{subsection}

\end{section}

\begin{section}{A closer look at the stable partial order}
      We now study in more detail the relations \(B_{n,d}^\circ\), and their
    reflexive and transitive closures. We note that \(B_{2,d}\)
    is a chain. For \(n=3\) the relation \(B_{3,d}^\circ\), and the Hasse
    diagram for \(B_{3,d}\), looks as \prettyref{fig:HasseS3d}.
    
  \begin{figure}[tbh]
    \begin{center}
      \setlength{\unitlength}{0.8cm}
      \subfigure[The graph of \(B_{3,d}^\circ\).
        Multiplication
        with \(\frac{x_1}{x_3}\) is 
        \(\rightarrow\), multiplication with \(\frac{x_2}{x_3}\) is
        \(\uparrow\), 
        multiplication with \(\frac{x_1}{x_2}\) is \(\searrow\).
      ]{
          \begin{picture}(10,9)    \thicklines
    \put(0,0){\vector(1,0){8}}
    \put(0,0){\vector(0,1){8}}
    \put(1,0){\circle*{0.3}}
    \put(0,1){\circle*{0.3}}
    \put(2,0){\circle*{0.3}}
    \put(1,1){\circle*{0.3}}
    \put(0,2){\circle*{0.3}}
    \put(3,0){\circle*{0.3}}
    \put(2,1){\circle*{0.3}}
    \put(1,2){\circle*{0.3}}
    \put(0,3){\circle*{0.3}}
    \put(4,0){\circle*{0.3}}
    \put(3,1){\circle*{0.3}}
    \put(2,2){\circle*{0.3}}
    \put(1,3){\circle*{0.3}}
    \put(0,4){\circle*{0.3}}
    \put(5,0){\circle*{0.3}}
    \put(4,1){\circle*{0.3}}
    \put(3,2){\circle*{0.3}}
    \put(2,3){\circle*{0.3}}
    \put(1,4){\circle*{0.3}}
    \put(0,5){\circle*{0.3}}
    \put(6,0){\circle*{0.3}}
    \put(5,1){\circle*{0.3}}
    \put(4,2){\circle*{0.3}}
    \put(3,3){\circle*{0.3}}
    \put(2,4){\circle*{0.3}}
    \put(1,5){\circle*{0.3}}
    \put(0,6){\circle*{0.3}}
    \put(7,0){\circle*{0.3}}
    \put(6,1){\circle*{0.3}}
    \put(5,2){\circle*{0.3}}
    \put(4,3){\circle*{0.3}}
    \put(3,4){\circle*{0.3}}
    \put(2,5){\circle*{0.3}}
    \put(1,6){\circle*{0.3}}
    \put(0,7){\circle*{0.3}}
    \put(0,0){\vector(1,0){0.9}}
    \put(0,0){\vector(0,1){0.9}}
    \put(1,0){\vector(1,0){0.9}}
    \put(1,0){\vector(0,1){0.9}}
    \put(0,1){\vector(1,0){0.9}}
    \put(0,1){\vector(0,1){0.9}}
    \put(2,0){\vector(1,0){0.9}}
    \put(2,0){\vector(0,1){0.9}}    
    \put(1,1){\vector(1,0){0.9}}
    \put(1,1){\vector(0,1){0.9}}
    \put(0,2){\vector(1,0){0.9}}
    \put(0,2){\vector(0,1){0.9}}
    \put(3,0){\vector(1,0){0.9}}
    \put(3,0){\vector(0,1){0.9}}    
    \put(2,1){\vector(1,0){0.9}}
    \put(2,1){\vector(0,1){0.9}}    
    \put(1,2){\vector(1,0){0.9}}
    \put(1,2){\vector(0,1){0.9}}    
    \put(0,3){\vector(1,0){0.9}}
    \put(0,3){\vector(0,1){0.9}}    
    \put(4,0){\vector(1,0){0.9}}
    \put(4,0){\vector(0,1){0.9}}    
    \put(3,1){\vector(1,0){0.9}}
    \put(3,1){\vector(0,1){0.9}}    
    \put(2,2){\vector(1,0){0.9}}
    \put(2,2){\vector(0,1){0.9}}    
    \put(1,3){\vector(1,0){0.9}}
    \put(1,3){\vector(0,1){0.9}}    
    \put(0,4){\vector(1,0){0.9}}
    \put(0,4){\vector(0,1){0.9}}    
    \put(5,0){\vector(1,0){0.9}}
    \put(5,0){\vector(0,1){0.9}}    
    \put(4,1){\vector(1,0){0.9}}
    \put(4,1){\vector(0,1){0.9}}    
    \put(3,2){\vector(1,0){0.9}}
    \put(3,2){\vector(0,1){0.9}}    
    \put(2,3){\vector(1,0){0.9}}
    \put(2,3){\vector(0,1){0.9}}    
    \put(1,4){\vector(1,0){0.9}}
    \put(1,4){\vector(0,1){0.9}}    
    \put(0,5){\vector(1,0){0.9}}
    \put(0,5){\vector(0,1){0.9}}    
    \put(6,0){\vector(1,0){0.9}}
    \put(6,0){\vector(0,1){0.9}}    
    \put(5,1){\vector(1,0){0.9}}
    \put(5,1){\vector(0,1){0.9}}    
    \put(4,2){\vector(1,0){0.9}}
    \put(4,2){\vector(0,1){0.9}}    
    \put(3,3){\vector(1,0){0.9}}
    \put(3,3){\vector(0,1){0.9}}    
    \put(2,4){\vector(1,0){0.9}}
    \put(2,4){\vector(0,1){0.9}}    
    \put(1,5){\vector(1,0){0.9}}
    \put(1,5){\vector(0,1){0.9}}    
    \put(0,6){\vector(1,0){0.9}}
    \put(0,6){\vector(0,1){0.9}}    
    \put(0,7){\vector(1,-1){0.9}}
    \put(1,6){\vector(1,-1){0.9}}
    \put(2,5){\vector(1,-1){0.9}}
    \put(3,4){\vector(1,-1){0.9}}
    \put(4,3){\vector(1,-1){0.9}}
    \put(5,2){\vector(1,-1){0.9}}
    \put(6,1){\vector(1,-1){0.9}}
  \end{picture}
        }      
      \setlength{\unitlength}{0.4cm}
      \subfigure[The Hasse diagram of \(B_{3,d}\).
      ]{
        \begin{picture}(8,15) \thicklines
\put(8,8){\line(-1,1){8}}
\put(0,0){\line(1,1){8}}
\put(0,2){\line(1,1){7}}
\put(0,4){\line(1,1){6}}
\put(0,6){\line(1,1){5}}
\put(0,8){\line(1,1){4}}
\put(0,10){\line(1,1){3}}
\put(0,12){\line(1,1){2}}
\put(0,14){\line(1,1){1}}
\put(0,0){\line(0,1){14}}
\put(1,1){\line(0,1){12}}
\put(2,2){\line(0,1){10}}
\put(3,3){\line(0,1){8}}
\put(4,4){\line(0,1){6}}
\put(5,5){\line(0,1){4}}
\put(6,6){\line(0,1){2}}
\put(0,0){\circle*{0.6}}
\put(0,2){\circle*{0.6}}
\put(0,4){\circle*{0.6}}
\put(0,6){\circle*{0.6}}
\put(0,8){\circle*{0.6}}
\put(0,10){\circle*{0.6}}
\put(0,12){\circle*{0.6}}
\put(0,14){\circle*{0.6}}
\put(0,16){\circle*{0.6}}
\put(1,1){\circle*{0.6}}
\put(1,3){\circle*{0.6}}
\put(1,5){\circle*{0.6}}
\put(1,7){\circle*{0.6}}
\put(1,9){\circle*{0.6}}
\put(1,11){\circle*{0.6}}
\put(1,13){\circle*{0.6}}
\put(1,15){\circle*{0.6}}
\put(2,2){\circle*{0.6}}
\put(2,4){\circle*{0.6}}
\put(2,6){\circle*{0.6}}
\put(2,8){\circle*{0.6}}
\put(2,10){\circle*{0.6}}
\put(2,12){\circle*{0.6}}
\put(2,14){\circle*{0.6}}
\put(3,3){\circle*{0.6}}
\put(3,5){\circle*{0.6}}
\put(3,7){\circle*{0.6}}
\put(3,9){\circle*{0.6}}
\put(3,11){\circle*{0.6}}
\put(3,13){\circle*{0.6}}
\put(4,4){\circle*{0.6}}
\put(4,6){\circle*{0.6}}
\put(4,8){\circle*{0.6}}
\put(4,10){\circle*{0.6}}
\put(4,12){\circle*{0.6}}
\put(5,5){\circle*{0.6}}
\put(5,7){\circle*{0.6}}
\put(5,9){\circle*{0.6}}
\put(5,11){\circle*{0.6}}
\put(6,6){\circle*{0.6}}
\put(6,8){\circle*{0.6}}
\put(6,10){\circle*{0.6}}
\put(7,7){\circle*{0.6}}
\put(7,9){\circle*{0.6}}
\put(8,8){\circle*{0.6}}

\put(0,14){\line(1,1){1.5}}
\end{picture}
        }
    \end{center}
      \caption{The stable relation and partial order, for \(n=3\).}
      \label{fig:HasseS3d}
    \end{figure}
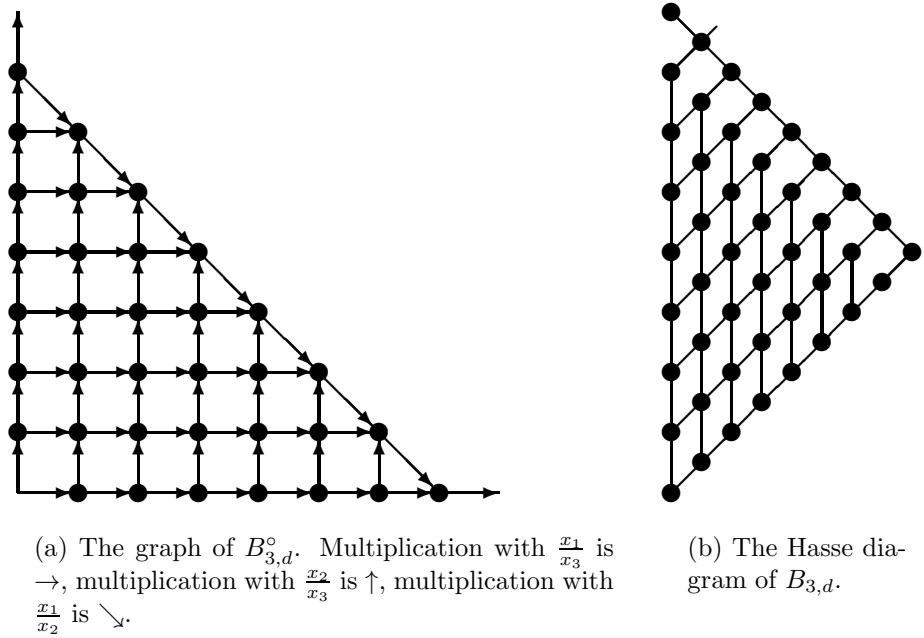

  
  
  \begin{definition}\label{def:End}
    We define
  \begin{displaymath}
    E_{n-1,d} = \setsuchas{x_1^{\alpha_1}\cdots
    x_{n-1}^{\alpha_{n-1}}}{\exists \alpha_n: x_1^{\alpha_1}\cdots
    x_{n-1}^{\alpha_{n-1}} x_n^{\alpha_n} \in \M_d^n} =
  \cup_{v=0}^d \M_v^{n-1}.
  \end{displaymath}
  \end{definition}
  As in \prettyref{fig:HasseS3d}, we can draw the graph of \(B_{n,d}\)
  as a graph of the vertex set \(E_{n-1,d}\). From this picture it
  can be seen that

  \begin{lemma}\label{lemma:Epart}
    \(M_d^n\) is the disjoint union of \(\M_d^{n-1}\) and \(x_n
    \M_{d-1}^n\). We have that
    \begin{align*}
      (\M_d^{n-1}, \restRel{B_{n,d}}{\M_d^{n-1}}) & \simeq
      (\M_d^{n-1}, B_{n-1,d}) \\
      (x_n \M_{d-1}^n, \restRel{B_{n,d}}{x_n \M_{d-1}^n}) &
      \simeq 
      (E_{n-1,d}, \restRel{D}{E_{n-1,d}})
    \end{align*}
    If \(m \in \M_d^{n-1}\), \(m' \in x_n \M_{d-1}^n\) then 
    \((m',m) \not \in B_{n,d}\), and \((m,m') \in B_{n,d}^\circ\)
    iff \(m = (x_i/x_n) m'\) for some \(1 \le i < n\). 
  \end{lemma}
  \begin{proof}
    We have that 
    \[(\M_d^{n-1},  \restRel{B_{n,d}}{\M_d^{n-1}})  
    \simeq  (\M_d^{n-1}, B_{n-1,d}).\]
    If \((p,p') \in B_{n,d}\) and \(p'\) is
    divisible by \(x_n\), then the relation between \(p\) and \(p'\)
    must be a substitution \(x_n \mapsto x_i\); if \(p \neq p'\) then
    \(1 \le i \le n-1\). This proves the last
    assertions.   
  \end{proof}

  \begin{theorem}\label{thm:Blattice}
    For any positive integers \(v,d\), the poset
    \((\M_d,B_{\cdot,d})\) is a lattice, and the poset 
    \((\M_d^v,B_{v,d})\) is a sublattice.
  \end{theorem}
  \begin{proof}
    It is clearly enough to show that 
    \((\M_d^n,B_{n,d})\) is a lattice for all \(n\). We do this
    by induction on \(n\), the case \(n=2\) is proved by the
    observation above that \((\M_d^2,B_{2,d})\) is a
    chain. Suppose that 
    \((\M_d^w,B_{w,d})\) is a lattice for 
    \(w < n\), then clearly 
    \(B_{i,d}\) is a sublattice of \(B_{j,d}\) for \(i \le j <
    n\). 
    Recall that a finite poset with a unique minimal and a unique maximal
    element is a lattice iff every pair of elements in \(\M_d^n\) has an
    infimum \cite{GLT}. So, to show that \((\M_d^n,B_{n,d})\) is a
    lattice, we    
    take \(m,m' \in \M_d^n\),  and want to define the infimum 
    \(m \wedge m' \in \M_d^n\). The following  cases present
    themselves: 
    \begin{enumerate}[(A)]
    \item \label{enum:case1}
      \(m,m' \in \M_d^{n-1}\). Then we use the induction hypothesis to
      define \(m \wedge m' \in \M_d^{n-1} \subset \M_d^n\).
    \item \label{enum:case2}
      \(m,m' \in \M_d^n \setminus \M_d^{n-1} = x_n \M_{d-1}^n\).
      Then we define the bijective map 
      \begin{align*}
      f: \M_d^n &\to E_{n-1,d}  \\
      x_i & \mapsto 
      \begin{cases}
        x_i & \text{ if } 1 \le i < n \\
        1 & \text{ if } i = n
      \end{cases}
      \end{align*}
      extended multiplicatively to all elements in \(\M_d^n\).
      We define
      \[m \wedge m' = f^{-1}(\gcd(f(m) , f(m'))),\]
      the gcd denotes the ordinary
      greatest common divisor on \(\M^{n-1}\). From
      \prettyref{lemma:Epart} it follows that this
      is indeed the
      greatest lower bound of \(m\) and \(m'\) in
      \[(x_n\M_{d-1}^{n},\restRel{B_{n,d}}{x_n\M_{d-1}^{n}}) 
      \simeq (E_{n-1,d}, \restRel{D}{E_{n-1,d}}).\] 
    \item \label{enum:case3}      
      \(m \in \M_d^{n-1}\), \(m' \in \M_d^n \setminus \M_d^{n-1}\).
      In this case, 
      first note that it is impossible that \(m'\) is larger than
      \(m\), because there is no transformation in \(B_{n,d}\) which
      introduces a \(x_n\). Furthermore, if \(v=\maxsupp(m)<n\) then
      the element \(m''=\frac{x_n}{x_v}m\) is the largest element in
      \(\M_d^n \setminus \M_d^{n-1}\) which is smaller than
      \(m\). Hence, the infimum of \(m''\) and \(m'\) is also the
      infimum of \(m\) and \(m'\).
    \end{enumerate}
  \end{proof}
  
  Recall \cite{Birkhoff:LT} that a lattice is non-modular iff it
  contains a copy of the non-modular lattice \(N_5\).
  \begin{figure}[htbp]
    \begin{center}
      \setlength{\unitlength}{0.4cm}
            \begin{picture}(2,4)
        \put(1,0){\circle*{0.3}}
        \put(0,1){\circle*{0.3}}
        \put(2,2){\circle*{0.3}}
        \put(0,3){\circle*{0.3}}
        \put(1,4){\circle*{0.3}}
        \put(1,0){\line(-1,1){1}}
        \put(1,0){\line(1,2){1}}
        \put(0,1){\line(0,1){2}}
        \put(0,3){\line(1,1){1}}
        \put(2,2){\line(-1,2){1}}
      \end{picture}
      \caption{\(N_5\)}
      \label{fig:N5}
    \end{center}
  \end{figure}
  
  \begin{theorem}
    The lattice \((\M^n_d, B_{n,d})\) is non-modular for
    \(n \ge 3\), \(d \ge 2\).
  \end{theorem}
  \begin{proof}
    \((\M^n_d, B_{n,d})\)
    contains \((\M_d^3, B_{3,d})\) which is non-modular for \(d
    \ge 2\), since 
    its Hasse diagram
    contains a copy of \(N_5\), as
    \prettyref{fig:HasseS3d} shows. 
  \end{proof}

\end{section}

\begin{section}{Filters in the strongly stable poset}
      \begin{subsection}{On the number of Borel sets}
    We shall calculate the number of filters in \(A_{n,d}\) with 
    a fixed cardinality \(v\). Our motivation is the article by Marinari and
    Ramella \cite{Marinari:SomeProp}, where the uniqueness of Borel subsets of 
    cardinality \(v\) of \(A_{3,d}\), for certain \(v\), is used to determine
    the possible 
    numerical resolutions of Borel ideals. 

    We start with the following simple
    observation: 
    \begin{theorem}\label{thm:Filtrec}
      Let \((T, \ge)\) be a finite poset. For any \(A \subset T\) and
      \(v \in \Nat\), we denote by \(\Phi(A,v)\) the number of filters
      of cardinality \(v\)
      in the sub-poset \((A, \ge)\). For any \(x \in A\) we put \(I(x)
      = \setsuchas{y \in A}{y \le x}\) and \(F(x) = \setsuchas{y \in
        A}{y \ge x}\). Then the following recurrence relation holds:
      \begin{equation}
        \label{eqn:reccard}
        \forall x \in A: \quad 
        \Phi(A,v) = \Phi(A \setminus I(x), v) + \Phi(A \setminus F(x),
        v - \#F(x)).
      \end{equation}
      Furthermore \(\Phi(A,0) =   \Phi(A,\#A)=1\).

      Denote by \(\phi(A)\) the number of filters in \((A,\ge)\), so
      that 
      \begin{math}
        \phi(A) = \sum_{v=0}^{\#A} \Phi(A,v).
      \end{math}
      Then 
      \begin{equation}
        \label{eqn:rec}
        \forall x \in A: \quad 
        \phi(A) = \phi(A \setminus I(x)) + \phi(A \setminus F(x)).
      \end{equation}
    \end{theorem}
    \begin{proof}
      It suffices to prove \prettyref{eqn:reccard}.
      Let \(P\) be any filter in \((A,\ge)\) of cardinality \(v\). If
      \(x \not \in P\) then 
      \(I(x) \cap P = \emptyset\), hence \(P\) is a filter in 
      \((A \setminus I(x), \ge)\), with cardinality \(v\). There are 
      \(\Phi(A \setminus I(x), v)\) such filters.

      If on the other hand \(x \in P\) then \(F(x) \subset
      P\). Clearly \(P \setminus F(x)\) is a filter of cardinality 
      \(v  - \#F(x)\) in \(A \setminus F(x)\), and conversely, if
      \(P'\) is such a filter, then \(P' \cup F(x)\) is a filter of
      cardinality \(v\) in \(A\). Thus there is a bijection between
      filters in \(A\) containing \(x\) and having cardinality \(v\),
      and filters in \(A \setminus F(x)\) of cardinality 
      \(v - \#F(x)\). There are 
      \(\Phi(A \setminus F(x), v - \#F(x))\) such filters.
    \end{proof}

    We now apply this to the poset \((\M_d^n,A_{n,d})\). For
    \(n=2\), this is a chain with \(d+1\) elements, hence 
    \begin{align*}
      \Phi(\M_d^2,v) & = 
      \begin{cases}
        1 & \text{ if } 0 \le v \le d+1\\
        0 & \text{ otherwise.}
      \end{cases}\\
      \phi(\M_d^2) & = d+2
    \end{align*}
    For \(n=3\), we get
      \begin{corr}\label{corr:filtn2}
        Denote by \(F(n,d,v)\) the number of filters of cardinality
        \(v\) in \((\M_d^n, A_{n,d})\), and by \(f(n,d)\) the
        number of filters in \((\M_d^n,A_{n,d})\). Then 
        \begin{align*}
          F(3,d,v) & = F(3,d-1,v) + F(3,d-1,v - (d+1)) \\
          F(3,1,v) & = 
          \begin{cases}
            1 & \text{ if } 0 \le v \le 3 \\
            0 & \text{ otherwise}
          \end{cases}\\
          f(3,d) & = 2 f(3,d-1) \\
          f(3,1) & = 4
        \end{align*}
        In fact, \(f(3,d) = 2^{d+1}\).
      \end{corr}
      \begin{proof}
        We have that
    \begin{align*}
        \M_d^3 \setminus I(x_2^d) &= x_1 \M_{d-1}^3 \\
        F(x_2^d) &= \M_d^2 \\
        \#F(x_2^d) &= d+1 \\
        \M_d^3 \setminus F(x_2^d) &= x_3 \M_{d-1}^3 
      \end{align*}
      We observe that there is a bijection between the filters in 
      \(\M_{d-1}^n\) of cardinality \(w\) and the filters in  \(t \M_{d-1}^n\) 
      of cardinality \(w\),
      for any  
      \(t \in  \M_1^n\); this bijection is given simply by multiplcation by
      \(t\).

      Hence, from \prettyref{thm:Filtrec} we get that
      \begin{align*}
        F(3,d,v) & =\Phi(\M_d^3,v) \\
        &= \Phi(\M_d^3 \setminus I(x_2^d),v) +
        \Phi(\M_d^3 \setminus F(x_2^d), v - \#F(x_2^d) \\
        & = \Phi(x_1 \M_{d-1}^3,v) + \Phi(x_3 \M_{d-1}^3, v-d-1) \\
        &= \Phi(\M_{d-1}^3,v) + \Phi(\M_{d-1}^3,v-d-1) \\
        &= F(3,d-1,v) + F(3,d-1,v-d-1)
      \end{align*}
      This establishes the recurrence relation. Since \(\M_1^3\) is a
      3-element chain, it has exactly one chain of length \(v\) when \(v \in
      \set{0,1,2,3}\), and no longer chains. Hence, the boundary conditions
      are as stated in the theorem.
      \end{proof}

    \begin{subsubsection}{``Splicing the filters''}
      We can give
      a recursion formula that is more
      efficient than \prettyref{eqn:reccard}. To elaborate on this
      formula, let for the remainder of this subsection \(n \ge 3\).  

      \begin{definition}\label{def:shrink}
        For \(F \subset \M^{n-1}\), define the \emph{interior} of
        \(F\) as 
        \begin{displaymath}
        \shrink(F) = \setsuchas{v \in F}{\forall i<j:
          \divides{x_i}{v} \implies   v (x_j/x_i) \in F},
        \end{displaymath}
        and the \emph{boundary}   of \(F\) as
      \begin{displaymath}
        \partial F = \setsuchas{m \in F}{\exists t \in \M^{n-1}
          \setminus 
          F: \exists i < j: m = (x_i/x_j) t}.
      \end{displaymath}
      \end{definition}
      
      \begin{theorem}\label{thm:skrinkfilt}
        Let \(F \subset \M_d^n\), define \(F_{d+1} =
        \emptyset\), and          
        \begin{displaymath}
        F_i := \setsuchas{m \in \M_{d-i}^{n-1}}{x_n^i m \in F}
        \text{ for }
        0 \le i \le d, 
        \end{displaymath}
        Then \(F\) is a filter in \((\M_d^n, A_{n,d})\) iff the
        following two conditions hold,  
        for \(0 \le i \le d\):
      \begin{enumerate}
      \item  \(F_i\)
        is a filter in \((\M_{d-i}^{n-1}, A_{n-1,d-i})\), 
      \item \(\shrink(F_i) \supset x_1 F_{i+1}\).
      \end{enumerate}
      \end{theorem}

      \begin{proof}
        Suppose first that \(F\) is a filter in \(\M_d^n\). 
        Let \(0 \le i \le d\). We prove that \(F_i\) is a filter in
        \(\M_{d-i}^{n-1}\). Namely, take \(a \in F_i\), so that
        \(x_n^i a \in F\), and let \(b \ge a\) (\wrt the
        strongly stable partial order) with \(b \in
        \M_{d-i}^{n-1}\). Then \(x_n^i b \in \M_d^n\), hence
        \(x_n^i b \ge x_n^i a\). Since \(F\) 
        is a filter, it follows that \(x_n^i b \in F\), hence that \(b 
        \in F_i\).

        We must also show that \(\shrink(F_i) \supset x_1
        F_{i+1}\). This condition is trivially fulfilled for \(i=d\),
        so suppose that \(1 \le i < d\). Take \(a = x_1 b \in x_1
        F_{i+1}\), so that \(x_n^{i+1} b \in F\). 
        To show that \(a \in \shrink(F_i)\), take \(r > s\) such that
        \(\divides{x_s}{b}\). We must 
        prove that \(\frac{x_r}{x_s}a \in F_i\), in other words, that
        \begin{equation}
          \label{eqn:crunx}
        \frac{x_r}{x_s} a x_n^i = \frac{x_r}{x_s} x_1 x_n^i b \in F.
        \end{equation}
        Clearly, there is a chain of elementary moves
        \begin{displaymath}
          x_s x_n \rightarrow x_s x_r \rightarrow x_1 x_r.
        \end{displaymath}
        Hence \(x_1 x_r \ge x_s x_n\), and hence 
        \(\frac{x_r}{x_s} x_1 \ge x_n\) in the ordered difference   group.
        Therefore, \(\frac{x_r}{x_s} x_1 x_n^i b \ge
          x_n^{i+1} b\), and since \(F\) is a filter, and since we
        have established that \(x_n^{i+1} b \in F\),
        \prettyref{eqn:crunx} follows. 
        The necessity of the conditions 1 and 2 is established.

        To prove sufficiency, suppose that \(F \subset \M_d^n\), and
        that conditions 1 and 2 hold. Suppose furthermore that \(t \in 
        F\), and that \(y \ge t\). Without loss of generality, we can
        assume that \(y\) is obtained from \(t\) by a single
        elementary move, so that \(t = z x_j\), \(y = z x_i\), \(i <
        j\). We distinguish between two cases: \(j < n\) and \(j =
        n\).

        If \(j < n\) we write \(z = z' x_n^{\ell}\), where
        \(\dividesnot{x_n}{z'}\). Since \(x_j z \in F\) it follows
        that \(x_j  z' \in F_{\ell}\). Using the fact that \(F_{\ell}\) is a
        filter (condition 1), we get that \(x_i z' \in F_{\ell}\), hence
        that \(y = x_i z 
        \in F\).

        There remains the case when \(j = n\). We write \(t = x_n z\), 
        \(y = x_i z\), \(i < n\). We also write \(z = z' x_n^{\ell}\) with
        \(\dividesnot{x_n}{z'}\). Then 
        \(t = x_n z = x_n^{{\ell}+1} z'\in F\),  
        hence \(z' \in F_{{\ell}+1}\). From the assumptions (condition 2)
        we get that 
        \begin{displaymath}
          x_1 z' \in x_1 F_{{\ell}+1} \subset \shrink(F_{\ell}),
        \end{displaymath}
        hence 
        \begin{displaymath}
          \frac{x_r}{x_s} x_1 z' \in F_{\ell} \quad \text{ for all } r > s \text{ 
            and } \divides{x_s}{x_1 z'}.
        \end{displaymath}
        In particular, \(\frac{x_i}{x_1} x_1 z' = x_i z' \in F_{\ell}\), hence 
        \begin{displaymath}
          F \ni x_n^{\ell} x_i z' = x_n^{\ell} z' x_i  = z x_i = y.
        \end{displaymath}
      \end{proof}
\end{subsubsection}

\begin{note}
  Another characterisation is given in \cite{BigattiRobbiano97}. 
\end{note}


      We can give a
      precise interpretation of the numbers \(F(3,d,v)\). A very
      similar reasoning is used in  \cite{Marinari:SomeProp}.
      
      \begin{prop}\label{prop:distinct}
        \(F(3,d,v)\) is equal to the number of numerical partitions of 
        \(v\) into distinct parts not exceeding \(d+1\).
      \end{prop}
      \begin{proof}
        Specialising \prettyref{thm:skrinkfilt} to the case \(n=3\) we 
        get 
        that a subset \(S \subset \M_d^3\) is a filter iff the
        following two conditions hold:
        \begin{enumerate}
        \item \(\forall \, 0\le i \le d: \quad (1/x_3^i) (S \cap x_3^i
          \M_{d-i}^2)\)  
          is a filter in \(\M_{d-i}^2\), and
        \item \(\forall \, 0 \le i < d: \quad \#(S \cap x_3^i \M_{d-i}^2) 
          > \#(S \cap x_3^{i+1} \M_{d-i-1}^2)\).
        \end{enumerate}
        Clearly, \(\#(S \cap x_3^i \M_{d-i}^2) \le d + 1\), and 
        \begin{displaymath}
          \#S = \sum_{i=0}^d \#(S \cap x_3^i \M_{d-i}^2).           
        \end{displaymath}
        Hence, the filters in \(\M_d^3\) of cardinality \(v\) are in
        bijective correspondence with the number of partitions of
        \(v\) into distinct parts not exceeding \(d+1\).
      \end{proof}
      
      This correspondence is illustrated in Figure~\ref{fig:fpart1}. 
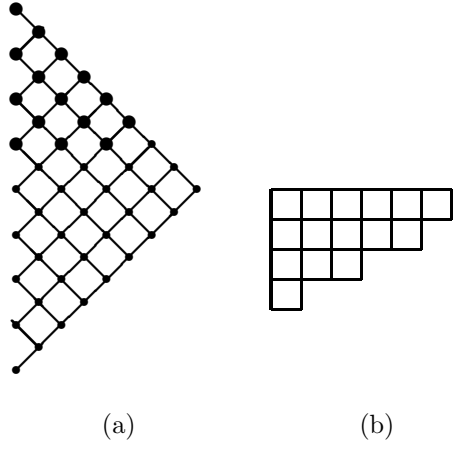
\begin{figure}[htbp]
  \begin{center}
      \setlength{\unitlength}{0.3cm}
      \subfigure[]{
        \begin{picture}(9,17) \thicklines
\put(1,1){\line(-1,1){1.2}}
\put(2,2){\line(-1,1){2}}
\put(3,3){\line(-1,1){3}}
\put(4,4){\line(-1,1){4}}
\put(5,5){\line(-1,1){5}}
\put(6,6){\line(-1,1){6}}
\put(7,7){\line(-1,1){7}}
\put(8,8){\line(-1,1){8}}
\put(0,0){\line(1,1){8}}
\put(0,2){\line(1,1){7}}
\put(0,4){\line(1,1){6}}
\put(0,6){\line(1,1){5}}
\put(0,8){\line(1,1){4}}
\put(0,10){\line(1,1){3}}
\put(0,12){\line(1,1){2}}
\put(0,14){\line(1,1){1.2}}
\put(0,0){\circle*{0.3}}
\put(0,2){\circle*{0.3}}
\put(0,4){\circle*{0.3}}
\put(0,6){\circle*{0.3}}
\put(0,8){\circle*{0.3}}
\put(0,10){\circle*{0.3}}
\put(0,12){\circle*{0.3}}
\put(0,14){\circle*{0.3}}
\put(0,16){\circle*{0.3}}
\put(1,1){\circle*{0.3}}
\put(1,3){\circle*{0.3}}
\put(1,5){\circle*{0.3}}
\put(1,7){\circle*{0.3}}
\put(1,9){\circle*{0.3}}
\put(1,11){\circle*{0.3}}
\put(1,13){\circle*{0.3}}
\put(1,15){\circle*{0.3}}
\put(2,2){\circle*{0.3}}
\put(2,4){\circle*{0.3}}
\put(2,6){\circle*{0.3}}
\put(2,8){\circle*{0.3}}
\put(2,10){\circle*{0.3}}
\put(2,12){\circle*{0.3}}
\put(2,14){\circle*{0.3}}
\put(3,3){\circle*{0.3}}
\put(3,5){\circle*{0.3}}
\put(3,7){\circle*{0.3}}
\put(3,9){\circle*{0.3}}
\put(3,11){\circle*{0.3}}
\put(3,13){\circle*{0.3}}
\put(4,4){\circle*{0.3}}
\put(4,6){\circle*{0.3}}
\put(4,8){\circle*{0.3}}
\put(4,10){\circle*{0.3}}
\put(4,12){\circle*{0.3}}
\put(5,5){\circle*{0.3}}
\put(5,7){\circle*{0.3}}
\put(5,9){\circle*{0.3}}
\put(5,11){\circle*{0.3}}
\put(6,6){\circle*{0.3}}
\put(6,8){\circle*{0.3}}
\put(6,10){\circle*{0.3}}
\put(7,7){\circle*{0.3}}
\put(7,9){\circle*{0.3}}
\put(8,8){\circle*{0.3}}
\multiput(0,16)(1,-1){6}{\circle*{0.6}}
\multiput(0,14)(1,-1){5}{\circle*{0.6}}
\multiput(0,12)(1,-1){3}{\circle*{0.6}}
\multiput(0,10)(1,-1){1}{\circle*{0.6}}
\end{picture}
      }
      \setlength{\unitlength}{0.4cm}
      \subfigure[]{
          \begin{picture}(7,5)    \thicklines
    \put(0,6){\line(1,0){6}}
    \multiput(0,6)(1,0){7}{\line(0,-1){1}}
    \put(0,5){\line(1,0){6}}
    \multiput(0,5)(1,0){6}{\line(0,-1){1}}
    \put(0,4){\line(1,0){5}}
    \multiput(0,4)(1,0){4}{\line(0,-1){1}}
    \put(0,3){\line(1,0){3}}
    \multiput(0,3)(1,0){2}{\line(0,-1){1}}
    \put(0,2){\line(1,0){1}}
  \end{picture}
      }
  \end{center}
\caption{A filter in \((\M_d^3, \rtr{A_{3,d}})\) and the associated
      numerical partition into distinct parts not exceeding \(d+1\).} 
  \label{fig:fpart1}
\end{figure}
\end{subsection}

\begin{subsection}{The poset structure of filters in the strongly
      stable partial order}
    \begin{subsubsection}{The case \(n=3\)}
      In fact, if we order the set of numerical partitions into
      distinct parts not exceeding \(d+1\) by inclusion of Young
      diagrams, and the filters in \(\M_d^3\) by inclusion, then the
      above bijection is easily seen to be a poset 
      isomorphism. There is yet another interpretation of this poset: 
      \begin{theorem}\label{thm:sqfree}
        The following three posets are isomorphic: 
        \begin{enumerate}[(i)]
        \item The filters in \((\M_d^3,A_{3,d})\), ordered by
          inclusion,
        \item The numerical partitions into distinct parts not
          exceeding \(d+1\), ordered by inclusion of Young diagrams,
        \item The poset \((\E^{d+1},
          \restRel{A_{d+1,\cdot}}{\E^{d+1}}))\), where 
          \(\E \subset \M\) denotes the set of all square-free
          monomials. 
        \end{enumerate}
        Furthermore, the following finite sets have the same cardinality:
        \begin{itemize}
        \item Filters in \((\M_d^3,A_{3,d})\) of cardinality
          \(v\),
        \item Numerical partitions of \(v\) into distinct parts not
          exceeding \(d+1\),
        \item Square-free monomials on \(\set{x_1,\dots,x_{d+1}}\) with
            weight \(v\), when \(x_i\) is given weight \(d+2-i\). 
        \end{itemize}
      \end{theorem}
      \begin{proof}
        Letting the variable
        \(x_i\) correspond to the singleton set \(\set{d+2-i}\), and a product 
        of distinct variables to the corresponding union, we regard a
        square-free monomial \(m \in \E^{d+1}\) as a 
        subset of \(\set{1,2,\dots,d+1}\). Summing
        the elements of this subset, we get a numerical partition into
        distinct parts not exceeding \(d+1\). This establishes a
        bijection, which we must show is isotone with isotone
        inverse. If \(m' = x_j m\) then the Young diagram of \(m'\)
        has an extra row compared to that of \(m\), and contains the
        latter; if \(m' = \frac{x_i}{x_{i+1}} m\) then one row of the Young
        diagram of \(m'\) is one unit longer than the corresponding row
        of the Young diagram of \(m\). The converse holds also, so the
        correspondence is an isomorphism.
      \end{proof}
      The Hasse diagram for the poset of square-free monomials in \(4\)
      variables is given in \prettyref{fig:sqfreePO}.
      \begin{figure}[htbp]
        \begin{center}
          \leavevmode
          \setlength{\unitlength}{0.8cm}
          \tiny{
            \begin{picture}(6,11) \thicklines
\put(3,0){\circle*{0.3}}
\put(3,0){\line(0,1){1}}
\put(3.3,0.2){\(1\)}

\put(3,1){\circle*{0.3}}
\put(3,1){\line(0,1){1}}
\put(3.3,1.1){\(x_4\)}

\put(3,2){\circle*{0.3}}
\put(3,2){\line(1,1){1}}
\put(3,2){\line(-1,1){1}}
\put(3.3,2.1){\(x_3\)}

\put(4,3){\circle*{0.3}}
\put(4,3){\line(-1,1){1}}
\put(4.3,3.1){\(x_3 x_4\)}

\put(2,3){\circle*{0.3}}
\put(2,3){\line(-1,1){1}}
\put(2,3){\line(1,1){1}}
\put(1.3,3){\(x_2\)}

\put(1,4){\circle*{0.3}}
\put(1,4){\line(1,1){1}}
\put(0.3,4){\(x_1\)}

\put(3,4){\circle*{0.3}}
\put(3,4){\line(-1,1){1}}
\put(3,4){\line(1,1){1}}
\put(3.3,4.1){\(x_2 x_4\)}

\put(2,5){\circle*{0.3}}
\put(2,5){\line(1,1){1}}
\put(0.7,5){\(x_1x_4\)}

\put(4,5){\circle*{0.3}}
\put(4,5){\line(-1,1){1}}
\put(4,5){\line(1,1){1}}
\put(4.3,5.1){\(x_2 x_3\)}

\put(3,6){\circle*{0.3}}
\put(3,6){\line(-1,1){1}}
\put(3,6){\line(1,1){1}}
\put(1.7,6){\(x_1x_3\)}

\put(5,6){\circle*{0.3}}
\put(5,6){\line(-1,1){1}}
\put(5.3,6.1){\(x_2 x_3 x_4\)}

\put(4,7){\circle*{0.3}}
\put(4,7){\line(-1,1){1}}
\put(4.3,7.1){\(x_1 x_3 x_4\)}

\put(2,7){\circle*{0.3}}
\put(2,7){\line(1,1){1}}
\put(0.7,7){\(x_1x_2\)}

\put(3,8){\circle*{0.3}}
\put(3,8){\line(0,1){1}}
\put(3.3,8.1){\(x_1 x_2 x_4\)}

\put(3,9){\circle*{0.3}}
\put(3,9){\line(0,1){1}}
\put(3.3,9.1){\(x_1 x_2 x_3\)}

\put(3,10){\circle*{0.3}}
\put(3.3,10.2){\(x_1 x_2 x_3 x_4\)}
\end{picture}
          }
          \caption{The poset \((\E^4,\restRel{A_{4,\cdot}}{\E^4})\).}
          \label{fig:sqfreePO}
        \end{center}
      \end{figure}
      
      It is clear that we can extend
      \prettyref{thm:filters} as follows (with the same definition
      of strongly stable ideal)
      \begin{lemma}
        For any positive integer \(n\), the strongly stable monoid
        ideals in \(\E^n\) are precisely the filters in \((\E^n,
        \restRel{A_{n,\cdot}}{\E^n})\). 
      \end{lemma}

      It is proved in \cite{Aramova:Gotzman} that
      \prettyref{thm:galligo} has a counterpart in the exterior
      algebra; that is, generic initial ideals in the exterior algebra
      \(\wedge^n\) are strongly stable, hence their intersections with
      \(\E^n\) are filters in \((\E^n,\restRel{A_{n,\cdot}}{\E^n})\).
      We can state this as follows:
      \begin{theorem}
        The sets of monomials  of a generic initial ideals in the
        exterior algebra on \(n\) 
        variables (with coefficients in \(\C\)) is a filter in the
        poset of filters of the poset \(A_{3,n-1}\).
      \end{theorem}
    \end{subsubsection}

    \begin{subsubsection}{The case \(n > 3\)}
      It is straightforward  to use
      \prettyref{thm:skrinkfilt} to give a description of the set of
      filters in \((\M_d^n, A_{n,d})\) in terms of certain
      hyper-partitions (see \cite{MacMahon}). In particular, for
      \(n=4\) we get
      \begin{prop}
        The following two posets are isomorphic:
        \begin{enumerate}
        \item         
          The poset of filters in \((\M_d^4, A_{4,d})\), ordered
          by inclusion,
        \item
          The poset of planar partitions
          which are contained in a \((d+1) \times (d+1) \times (d+1)\)
          box, and which have horizontal and vertical ``steps'' of
          length \(\le 1\).
        \end{enumerate}
      \end{prop}
      An explanation of the nomenclature is in order: we draw the solid
      Young diagram as  a union of unit
      cubes \([i,i+1] \times [j,j+1] \times [k,k+1]\), with \(i,j,k\)
      non-negative integers satisfying \(0 \le k < \pi_{i,j}\); the
      \(\pi_{i,j}\)'s are non-negative integers, almost all zero, such
      that for all \(x,y\) we have that \(\pi_{x+1,y} \le \pi_{x,y}\)
      and \(\pi_{x,y+1} \le \pi_{x,y}\). Then we demand in addition
      that each ``level'' of cubes should represent a partition into
      distinct parts, and that each ``level'' should be contained in
      the interior of the one upon which it rests. It is easy to see
      that  this implies that there can be no
      vertical steps of height 2 or greater. Similarly, each
      ``level'', when drawn in this fashion has steps of height 1.
    \end{subsubsection}
  \end{subsection}

\end{section}

\begin{section}{Filters in the  stable poset}  
    \begin{subsection}{On the number of stable subsets of \(\M_d^3\)}
    Let us briefly consider the question of how many filters there are
    (of a given cardinality)
    in \((\M_d^n, B_{n,d})\). For \(n=1,2\) the partial order is 
    a chain, hence the enumeration of filters is trivial. For \(n=3\), 
    the Hasse diagram looks like \prettyref{fig:HasseS3d}. We apply
    \prettyref{thm:Filtrec}, by partitioning the filters into two
    classes: those that contain \(x_2^d\), and those that do not. We
    have that 
    \begin{align*}
      F(x_2^d) &= \M_d^2 \\
      \M_d^3 \setminus F(x_2^d) & = x_3 \M_{d-1}^3 \\
      \M_d^3 \setminus I(x_2^d) &= x_1 \M_{d-1}^3 \\
    \end{align*}
    It is easy to see that 
    \[(x_1 \M_{d-1}^3,\restRel{B_{3,d}}{x_1 \M_{d-1}^3})
    \simeq   (\M_{d-1}^3, B_{3,d-1}).\]
    Furthermore, in \(x_3 \M_{d-1}^3\) every monomial is divisible by
    \(x_3\), hence the allowed substitutions are \(m \mapsto x_1/x_3\) 
    and \(m \mapsto x_2/x_3\). Hence (recall \prettyref{def:End})
    \begin{displaymath}
      (x_3 \M_{d-1}^3, \restRel{B_{3,d}}{x_3 \M_{d-1}^3}) 
      \simeq (E_{2,d-1}, \restRel{D}{E_{2,d-1}})
    \end{displaymath}
    The poset \((E_{2,d},\restRel{D}{E_{2,d}})\) is a sub-poset of
    \(\Nat^2\)  
    with the natural divisibility order.

    So, if we denote by \(G(d)\) the number of filters in
    \((\M_d^3,B_{3,d})\), by \(GG(d,v)\) the number of such
    filters of cardinality \(v\), by \(C(d)\) the number of filters in
    \((E_{2,d},\restRel{D}{E_{2,d}})\), and by \(CC(d,v)\) the number
    of such  
    filters of cardinality \(v\), we have that
    \begin{align}\label{eqn:stabrecur}
      G(d) & = G(d-1) + C(d-1) \\
      GG(d,v) & = GG(d-1, v - d - 1) + CC(d-1,v)
    \end{align}
    
    We can give an illuminating interpretation of the numbers \(C(d)\)
    by observing that filters in \((E_{2,d},\restRel{D}{E_{2,d}})\)
    are in 
    bijective correspondence with lattice walks in \(E_{2,d+2}\) from
    \((0,d+2)\) to \((d+2,0)\), consisting of moves of unit length
    down and to the right. Namely, to such a walk we associate the
    filter 
    in \((E_d,\restRel{D}{E_{2,d}})\) which is generated by all lattice
    points in \(E_{2,d}\) which are visited during the walk. Thus, the
    empty filter corresponds to the walk
    \begin{displaymath}
      (0,d+2) \downarrow (0,d+1) \rightarrow (1,d+1) \downarrow
      \rightarrow \cdots \downarrow (d+1,0) \rightarrow (d+2,0)
    \end{displaymath}
    whereas the filter \(E_{2,d}\) corresponds to the walk
    \begin{displaymath}
      (0,d+2) \downarrow (0,d+1) \downarrow \cdots \downarrow (0,0)
      \rightarrow (1,0) \rightarrow \cdots \rightarrow (d+1,0)
      \rightarrow (d+2,0)
    \end{displaymath}
    The correspondence for a general filter is illustrated in
    \prettyref{fig:walkfilt}.
    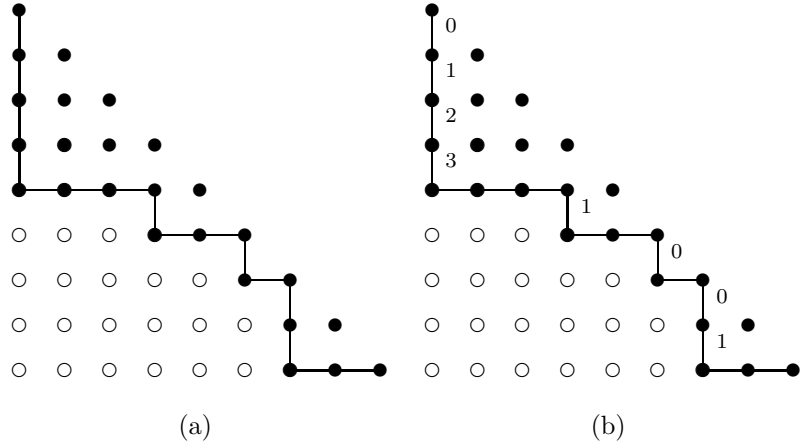
\begin{figure}[htbp]
      \begin{center}
        \subfigure[]{
      \setlength{\unitlength}{0.6cm}
      \begin{picture}(8,8)
\multiput(0,0)(0,1){7}{\circle{0.3}}
\multiput(1,0)(0,1){6}{\circle{0.3}}
\multiput(2,0)(0,1){5}{\circle{0.3}}
\multiput(3,0)(0,1){4}{\circle{0.3}}
\multiput(4,0)(0,1){3}{\circle{0.3}}
\multiput(5,0)(0,1){2}{\circle{0.3}}
\multiput(6,0)(0,1){1}{\circle{0.3}}

\put(7,0){\circle*{0.3}}
\put(6,1){\circle*{0.3}}
\put(5,2){\circle*{0.3}}
\put(4,3){\circle*{0.3}}
\put(3,4){\circle*{0.3}}
\put(2,5){\circle*{0.3}}
\put(1,6){\circle*{0.3}}
\put(0,7){\circle*{0.3}}

\put(8,0){\circle*{0.3}}
\put(7,1){\circle*{0.3}}
\put(6,2){\circle*{0.3}}
\put(5,3){\circle*{0.3}}
\put(4,4){\circle*{0.3}}
\put(3,5){\circle*{0.3}}
\put(2,6){\circle*{0.3}}
\put(1,7){\circle*{0.3}}
\put(0,8){\circle*{0.3}}

\put(0,8){\line(0,-1){4}}
\put(0,4){\line(1,0){3}}
\put(3,4){\line(0,-1){1}}
\put(3,3){\line(1,0){2}}
\put(5,3){\line(0,-1){1}}
\put(5,2){\line(1,0){1}}
\put(6,2){\line(0,-1){2}}
\put(6,0){\line(1,0){2}}

\put(0,4){\circle*{0.3}}
\put(0,5){\circle*{0.3}}
\put(0,6){\circle*{0.3}}
\put(1,4){\circle*{0.3}}
\put(1,5){\circle*{0.3}}
\put(2,4){\circle*{0.3}}
\put(3,3){\circle*{0.3}}
\put(6,0){\circle*{0.3}}
\end{picture}
      }
\subfigure[]{
      \setlength{\unitlength}{0.6cm}
      \begin{picture}(8,8)
\multiput(0,0)(0,1){7}{\circle{0.3}}
\multiput(1,0)(0,1){6}{\circle{0.3}}
\multiput(2,0)(0,1){5}{\circle{0.3}}
\multiput(3,0)(0,1){4}{\circle{0.3}}
\multiput(4,0)(0,1){3}{\circle{0.3}}
\multiput(5,0)(0,1){2}{\circle{0.3}}
\multiput(6,0)(0,1){1}{\circle{0.3}}

\put(7,0){\circle*{0.3}}
\put(6,1){\circle*{0.3}}
\put(5,2){\circle*{0.3}}
\put(4,3){\circle*{0.3}}
\put(3,4){\circle*{0.3}}
\put(2,5){\circle*{0.3}}
\put(1,6){\circle*{0.3}}
\put(0,7){\circle*{0.3}}

\put(8,0){\circle*{0.3}}
\put(7,1){\circle*{0.3}}
\put(6,2){\circle*{0.3}}
\put(5,3){\circle*{0.3}}
\put(4,4){\circle*{0.3}}
\put(3,5){\circle*{0.3}}
\put(2,6){\circle*{0.3}}
\put(1,7){\circle*{0.3}}
\put(0,8){\circle*{0.3}}

\put(0,8){\line(0,-1){4}} 
\put(0,4){\line(1,0){3}}
\put(3,4){\line(0,-1){1}}
\put(3,3){\line(1,0){2}}
\put(5,3){\line(0,-1){1}}
\put(5,2){\line(1,0){1}}
\put(6,2){\line(0,-1){2}}
\put(6,0){\line(1,0){2}}

\put(0,4){\circle*{0.3}}
\put(0,5){\circle*{0.3}}
\put(0,6){\circle*{0.3}}
\put(1,4){\circle*{0.3}}
\put(1,5){\circle*{0.3}}
\put(2,4){\circle*{0.3}}
\put(3,3){\circle*{0.3}}
\put(6,0){\circle*{0.3}}
\tiny{
\put(0.3,7.5){0}
\put(0.3,6.5){1}
\put(0.3,5.5){2}
\put(0.3,4.5){3}
\put(3.3,3.5){1}
\put(5.3,2.5){0}
\put(6.3,1.5){0}
\put(6.3,0.5){1}
}
\end{picture}
      }
      \caption{A lattice walk in \(E_{2,8}\) and its associated filter
      in \(E_{2,6}\).} 
      \label{fig:walkfilt}
      \end{center}
    \end{figure}
    It is well-known that the number of lattice walks in \(E_{2,N+1}\)
    is the \(N\)'th Catalan number \(C_N = (2N)!/(N!(N+1)!)\). It
    follows that 
    \(C(d) = C_{d-1}\). Hence we can solve
    \prettyref{eqn:stabrecur} and get
    \begin{prop}
      The number of filters in \((\M_d^3,B_{3,d})\) is 
      \begin{math}
        \sum_{i=0}^{d+1} C_i,
      \end{math}
      where \(C_i\) is the \(i\)'th Catalan number.
    \end{prop}

    The following lemma shows that we can find an explicit, although
    complicated, recurrence relation for \(CC(d,v)\), and hence for
    \(GG(d,v)\). 
    \begin{lemma}
      We label the horizontal edges connecting vertices in
      \(E_{2,d+2}\) with 0, and the vertical edges \((a,b) \to
      (a,b-1)\) with \(d+2-a-b\). For a path using only such edges, we
      define its weight to be the sum of the labels of the edges.
      Let \(S_d(a,b,w)\) be the number of such paths in \(E_{2,d+2}\)
      from \((a,b)\) to  \((d+2,0)\) of weight \(w\).
      Then \(CC(d,v) =   S_d(0,d+2,w)\), and \(S_d(a,b,w)\) is the
      unique solution to the 
      following system of equations:
      \begin{equation}\label{eqn:wr}
        S_d(a,b,w) = 
        \begin{cases}
          0 & (a,b) \not \in E_{2,d+2} \\
          \begin{cases}
            1 & w = d - a + 1 \\
            0 & \text{ otherwise}
          \end{cases} & b=0 \\
          \begin{cases}
            1 & w \le d+1-a \\
            0 & \text{ otherwise}
          \end{cases} & b=1 \\
          \sum_{j=a}^{d+2-b} S_d(j,b-1, w+j+b-d-2) & \text{
          otherwise} 
        \end{cases}
      \end{equation}
    \end{lemma}
    \begin{proof}
      It is easy to see that the weight of the path counts the
      cardinality of the associated filter in \(E_{2,d}\). From
      \((a,b)\), one can choose to descend from 
      \begin{displaymath}
      (a,b),   (a+1,b),\dots,(d+2-b,b);        
      \end{displaymath}
      the corresponding vertical step
      \((j,b) \to (j,b-1)\) has weight \(d+2-b-j\). Hence,
      \(S_d(a,b,w) = \sum_{j=a}^{d+2-b} S_d(j,b-1, d+2-b-j)\). The
      boundary conditions are easily verified.
    \end{proof}
    
    We illustrate the weighting of the edges in
    \prettyref{fig:walkfilt}.
    
    In general, \prettyref{eqn:wr} seems hard to solve explicitly, but
    for the special case \(d > v\) we can indeed find the value of
    \(CC(d,v)\). 
    Recall \cite[Example 10.12]{asympt}
    that an \((n,k)\) fountain is an
    arrangement of \(n\) coins such that there are \(k\) coins in the
    bottom row, and such that each coin in a higher row rests on
    exactly two coins in the next lower row; a fountain of \(w\) coins
    is any \((w,k)\).
    \begin{prop}
      If \(d > w\) then the stable subsets of \(\M_d^3\) of cardinality
      \(w\) are in bijective correspondence with the fountains of
      \(w\) coins.
    \end{prop}
    \begin{proof}
      Since the principal filter on \(x_2^d\) contains \(d+1\)
      elements, it is clear that for \(w < d\) a filter \(F \subset
      \M_d^3\), \(\#F=w\)  \wrt the stable partial order
      never contains an element \(\le x_2^d\). Hence, such a filter
      can be identified with a filter in the inductive limit 
      \(G=\varinjlim_v (\M_v^3,A_{3,v})\), where the injections
      are given by 
      \begin{align*}
        i: \M_v^3 & \to \M_{v+1}^3 \\
        m & \mapsto x_1 m
      \end{align*}
      The infinite Hasse diagram of \(G\) looks like
      \prettyref{fig:HasseS3d}(b), only extended infinitely downwards and to
      the right. A finite
      filter in \(G\) gives a fountain of coins by reading the
      successive diagonal rows from right to left, and conversely.
    \end{proof}
    It is known \cite{asympt}(see Example 10.12) that if \(a_w\)
    denotes the number of \(w\)-fountains, then a generating
    function is given by
    \begin{displaymath}
      f(z) = \sum_{n=0}^\infty a_n z^n =
      \frac{1}{1-\displaystyle{\frac{z}{1-
            \displaystyle{\frac{z^2}{1-\displaystyle{\frac{z^3}{1\cdots}}}}}}} 
    \end{displaymath}
\end{subsection}

\end{section}

\bibliographystyle{plain}

\end{document}